\documentclass[a4paper,fleqn,final]{cas-sc}

\usepackage{algorithm,algpseudocode,subcaption,float,tikz,amsthm}
\usepackage{enumerate}
\usepackage[authoryear,longnamesfirst]{natbib}

\usetikzlibrary{
  decorations.pathreplacing,
  positioning,
  calc,
  intersections,
  3d,
  shapes.geometric,
  shapes,
  chains,
  math,
  fit,
  backgrounds
}

\newtheorem{remark}{Remark}
\newtheorem{proposition}{Proposition}
\newtheorem{lemma}{Lemma}

\newtheorem{theorem}{Theorem}
\newtheorem{definition}{Definition}

\newcommand\restr[2]{{%
    \left.\kern-\nulldelimiterspace
      #1
      \vphantom{\big|}
    \right|_{#2}
  }}

\definecolor{RED}{rgb}{1,0,0}
\newcommand{\PL}[1]{\textcolor{black}{#1}}

\renewcommand{\vec}[1]{\mathbf{#1}}
\newcommand{\Pideal}[0]{P_{\mathrm{ideal}}}
\newcommand{\Rideal}[0]{R_{\mathrm{ideal}}}
\newcommand{\Mideal}[0]{M_{0,\mathrm{ideal}}}
\newcommand{\Sideal}[0]{S^{-1}_{\mathrm{ideal}}}

\newcommand{\Ker}[0]{\mathrm{Ker}(E^{(c)})}

\begin{document}

\let\WriteBookmarks\relax

% === Running headers ===
\shorttitle{On Spectral Clustering in Algebraic Multigrid Methods}
\shortauthors{{Lucero Lorca} et al.}

% === Title ===
\title[mode=title]{On Spectral Clustering in Algebraic Multigrid Methods}

% === Authors ===
\author{José Pablo {Lucero Lorca}}

\author[1]{Conor McCoid}

\author[2]{Michal Outrata}
\fnmark[1]

% === Affiliations ===
\affiliation[1]{
  organization={McMaster University},
  state={Ontario},
  country={Canada}
}

\affiliation[2]{
  organization={Charles University},
  city={Prague},
  country={Czechia}
}

% === Funding acknowledgment ===
\fntext[1]{This work was supported by the grants UNCE/24/SCI/005 and
  PRIMUS/25/SCI/022 of Charles University.}

\begin{abstract}
  We introduce a new direct multilevel method for solving \PL{general
    nonsingular} complex square linear systems \PL{admitting a
    suitable block partition, using a} regular smoother and an
  arbitrary but equal number of pre- and post-smoothings. Through
  careful analysis of the error propagation operator, we cluster the
  spectrum of this operator. This allows us to write a direct K-cycle
  version of the method.
\end{abstract}


\maketitle

\section{Introduction}\label{sec_Intro}

Multigrid methods exploit hierarchy in discretization to approximate
solutions~\citep{Hackbusch1985, Wesseling1992}. These hierarchies are
constructed by combining smooth and coarse representations of the
problem on each level of the hierarchy so that at different levels the
method focuses on different components of the error. The success of
multigrid methods depends on the interaction of its components:
\begin{itemize}
\item the smoother;
\item the smoothing parameters;
\item restriction and prolongation operators, and
\item the coarse problems.
\end{itemize}
Despite their wide use, the processes inside multigrid methods are not
fully understood, especially for discrete non-Hermitian problems.  One
must often rely on case-dependent analysis or numerical
experimentation to obtain insight into the mechanism of the method for
specific problems \citep{LuceroLorcaGander2024,
  HemkerHoffmannVanRaalte2004}.

Often, the matrix of the discrete system to be solved with a multigrid
method is partitioned in a block structure with square diagonal
blocks.  If there are only two levels in the multigrid method, then
there exists a direct V-cycle method without post-smoothing and with a
singular smoother which is equivalent to the classic $2\times 2$-block
inversion formula \citep{MacLachlanManteuffelMcCormick2006,
  GarciaRamosKehlNabben2020, LuceroLorca2025,
  LuceroLorcaRosenbergJankovMcCoidGander2025}.  It remains open
whether one can always construct a direct two-level V-cycle with a
regular smoother.  Note that a two-level method can be generalized to
a multilevel method through recursion, and so resolving this problem
for two-level methods also resolves it for multilevel methods
\citep{MacLachlanManteuffelMcCormick2006, Notay2005, Vassilevski2008,
  TangNabbenVuikErlangga2009, GosslerNabben2016, LuceroLorca2025}.

The goal of this manuscript is to present a new direct multilevel
method for general \PL{nonsingular complex matrices admitting a block
  partition with nonsingular diagonal blocks}, not necessarily
\PL{arising} from a PDE discretization, featuring a regular
smoother. We do so by clustering the spectrum of a two-level algebraic
multigrid method that uses an arbitrary, but equal, number of pre- and
post-smoothings and optimal smoothing parameters for which we provide
a closed-form formula.  Relying on the results for Krylov subspace
methods, especially their analysis based on eigenproperties of the
system matrix, we then show that the K-cycle multigrid
method\footnote{Citing~\citep[p. 474]{NotayVassilevski2008}, we recall
  that ``\emph{with K-cycle, the MG method is still based on the
    recursive use of a two-grid method, but the needed coarse-grid
    solve is defined by a few steps of a Krylov subspace method of
    choice with the already defined (by recursion) MG method on the
    previous (coarser) level as the preconditioner.}''.} becomes
direct.  Following the same idea, we then use recursion to make the
method multilevel.  Running a K-cycle, we make a choice of the Krylov
subspace method, see~\citep{LiesenStrakos2013}. Classical choices
include the generalized minimal residual method (GMRES) or, if the
problem is Hermitian positive definite (HPD), the conjugate gradient
method (CG). Importantly, if a multilevel method uses a non-symmetric
cycle (i.e., using different number of pre and post smoothings), then
a system using such a method as a preconditioner can no longer use CG
in the HPD case, leading to a loss of efficiency.  Thus, in many
applications there is a need for a symmetric cycle.

To cluster the spectrum, we perform spectral analysis on the error
propagation operator of a general algebraic multigrid method. The
analysis requires investigating the generalized eigenvalues and
eigenvectors of the matrix pencil of the Schur complement.  We need to
assume the matrix pencil is diagonalizable so that $\mathbb{C}^N$ can
be decomposed into a direct sum of low-dimensional invariant subspaces
of the error propagation operator \PL{and a
  complement on which the coarse correction eliminates the error
  exactly}.  This allows us to treat each
invariant subspace separately, giving a closed formula for the
eigenvalues of the error propagation operator as a polynomial in the
smoothing parameters and the generalized eigenvalues of the matrix
pencil. Clustering of the eigenvalues of the error propagation
operator then becomes equivalent to removing the dependency of the
polynomial on the generalized eigenvalues -- a condition used to
select the smoothing parameters.

The derivation of this result is purely algebraic.  Therefore, it can
be applied regardless of the discretization.  The clustering of the
spectrum is a significant improvement over more standard fine-tuning
of smoothing parameters which is bespoke for each PDE.

The paper is structured as follows:
in Section \ref{sec:Background} we recall previous results to be used
here;
in Section \ref{sec:newclustered} we describe the choice of operators
and parameters that cluster the spectrum;
in Section \ref{sec:newdirect} we describe the new multilevel
method and its corresponding direct K-cycle;
in Section \ref{sec:spectralanalysis} we analyze the spectrum of the error
propagation operator for this new method;
in Section \ref{sec:spectralclustering} we derive the optimal choice of the
smoothing parameters that cluster the eigenvalues of this operator;
in Section \ref{sec:adaptation} we illustrate this clustering
with numerical evidence;
finally, in Section \ref{sec:conclusion} we summarize our novel results
and briefly describe future directions for this research.

\subsection{Existing methods}

Other methods that use hierarchies are domain decomposition methods
with coarse correction~\citep{DoleanJolivetNataf2015}, hierarchical
direct solvers~\citep{Hackbusch2015,Martinsson2015}, and many
parallel-in-time methods~\citep{Gander2015}.

Multigrid methods have a long history, beginning with analyses of
relaxation for elliptic difference equations \citep{Fedorenko1962,
  Fedorenko1964, Bakhvalov1966, Brandt1977, Hackbusch1985,
  Wesseling1992, BriggsHensonMcCormick2000}, then quickly followed by
their algebraic counterparts \citep{Brandt1986,
  BrandtMcCormickRuge1985, RugeStueben1987}.  Many authors analyze
algebraic multigrid methods only for the HPD setting and over the real
numbers, due to the convenient properties of the inner product induced
by the system matrix within the analysis \citep{Notay2005,
  Vassilevski2008, TangNabbenVuikErlangga2009,
  GosslerNabben2016}. Notably, in~\citep{GosslerNabben2016}, the
authors also pose the question of the optimal smoothing parameters but
their analysis differs fundamentally in the considered set-up and
hence also in their results.  In the HPD setting, the existence of an
overall optimal algorithm has been considered in
\citep{FalgoutVassilevski2004, XuZhang2018} by means of an
energy-minimizing interpolation.

The non-symmetric algebraic multigrid literature offers several
approaches for adapting two-level methods to challenging problems.
These include:
\begin{itemize}
\item approximating the ideal restriction and prolongation operators
  under fixed sparsity patterns with constrained optimization
  techniques \citep{OlsonSchroderTuminaro2011,
    WiesnerTuminaroWallGee2014};
\item Petrov--Galerkin smoothed aggregation \citep{SalaTuminaro2008,
    BrezinaManteuffelMcCormickRugeSanders2010};
\item AIR and nAIR algorithms \citep{ManteuffelRugeSouthworth2018};
\item using adaptive
  strategies~\citep{MacLachlanManteuffelMcCormick2006,
    BrannickFrommerKahlMacLachlanZikatanov2010}.
\end{itemize}
These approaches deliver iterative methods instead of direct ones, but
can potentially be used to adapt our new direct method to those same
problems.

\section{Background}\label{sec:Background}

For algebraic multigrid, the standard approach is to split the degrees
of freedom (DoFs) into two categories: fine and coarse.  This
splitting may be done by re-ordering the DoFs and separating the
matrix $L$ into blocks\footnote{In the AMG community the blocks are
  often labeled with subscripts $c$ (coarse) and $f$ (fine). Here, for
  ease of notation, we denote these blocks $A,\dotsc ,D$.}:
\begin{equation*}
  L = \begin{bmatrix}
        L_{ff} & L_{fc} \\
        L_{cf} & L_{cc}
      \end{bmatrix} \equiv
      \begin{bmatrix}
        A & B \\
        C & D
      \end{bmatrix},
  \quad L_{ff}\equiv A\in\mathbb{C}^{(N-n)\times (N-n)},
  \quad L_{cc}\equiv D\in\mathbb{C}^{n\times n}.
\end{equation*}
This partition of $L$ may arise from a physical decomposition, a
coloring or ordering of the DoFs or simply from algebraic convenience;
in all cases, we assume that the diagonal blocks represent subproblems
that can be inverted or approximately solved efficiently.  In such
case, the partitioning induces a closed formula for $L^{-1}$ due to
the Schur complement formulation of the LDU factorization of $L$.
These calculations can be written as a direct two-level algebraic
multigrid scheme and their approximation gives rise to a number of new
multigrid methods, see~\citep{MacLachlanManteuffelMcCormick2006,
  LuceroLorca2025, FoersterSteubenTrottenberg1981,
  RiesTrottenbergWinter1983, ChanWan2000}.

The general, symmetric two-level multigrid V-cycle preconditioner is
presented in Algorithm \ref{alg:TwoLevelVCycle}. This algorithm
includes $m$ pre- and post-smoothing steps, with each using the same
smoother $S^{-1}$ and smoothing parameters $\alpha_i \in \mathbb{C}$
symmetrically with respect to the coarse correction.  Also to be
defined are restriction, prolongation, and coarse space operators $R$,
$P$, and $M_0$, respectively. The operation of Algorithm
\ref{alg:TwoLevelVCycle} on $\mathbf{v} \in \mathbb{C}^N$ corresponds
to a linear operator, denoted by $M^{-1}$, acting on $\mathbf{v}$.

\begin{algorithm}
  \caption{TwoLevelVCycle($\mathbf{v}$)}
  \begin{algorithmic}[1]
    \State $\mathbf{x} \leftarrow \mathbf{0}$
    \For{$i=1,\dots,m$}
    \State $\mathbf{x} \leftarrow \mathbf{x} + \alpha_i S^{-1} (\mathbf{v} - L \mathbf{x})$
    \Comment{pre-smoothing}
    \EndFor
    \State $\mathbf{x} \leftarrow \mathbf{x} + P M_0^{-1} R (\mathbf{v} - L \mathbf{x})$
    \Comment{coarse correction}
    \For{$i=m,\dots,1$}
    \State $\mathbf{x} \leftarrow \mathbf{x} + \alpha_i S^{-1} (\mathbf{v} - L \mathbf{x})$
    \Comment{post-smoothing}
    \EndFor
    \State \textbf{return} \PL{$\mathbf{x}$}
  \end{algorithmic}
  \label{alg:TwoLevelVCycle}
\end{algorithm}

It is known (see \citep[Proposition 3.3 and Algorithm
1]{LuceroLorcaRosenbergJankovMcCoidGander2025} or~\citep[Corollary
2]{ManteuffelMuenzenmaierRugeSouthworth2019}) that Algorithm
\ref{alg:TwoLevelVCycle} without post-smoothing and with the choices
\begin{align}\label{eqn:optimalchoice1}
  \begin{aligned}
    P = \begin{bmatrix}
          -A^{-1} B \\
          I
        \end{bmatrix}, \quad
    R = \begin{bmatrix}
          \star & I
        \end{bmatrix}, \quad
    S^{-1} = \begin{bmatrix}
               A^{-1} \\
               &
             \end{bmatrix}, \quad
    m = 1, \quad
    \alpha_i = 1~\forall i,   \\
    \implies M_0 = R L P = D - C A^{-1} B,
  \end{aligned}
\end{align}
is a direct method, i.e., we have $\mathbf{x} = L^{-1}\mathbf{v}$
after executing line 5 of Algorithm~\ref{alg:TwoLevelVCycle}, for any
choice of $\star$. As line 7 is formulated as a smoothing correction,
we observe that even with post-smoothing
Algorithm~\ref{alg:TwoLevelVCycle} remains a direct method, i.e.,
with~\eqref{eqn:optimalchoice1} we have
$M^{-1} \mathbf{v} = L^{-1}\mathbf{v}$ for all vectors $\mathbf{v}$.

Algorithm \ref{alg:TwoLevelVCycle} is usually formulated in the
non-symmetric form, i.e., without the post-smoothing, for obvious
reasons -- the symmetrization on lines 6-8 is essentially
artificial. It is not obtained by necessity of the algebra but the
users' insistence on a symmetric cycle method (by the same token, one
often encounters $\star = -CA^{-1}$). Hence, when designing
preconditioners by approximation of direct solvers, the symmetrized
method lacks a natural target for the approximation. In many cases
this turns out to be of little concern as the symmetrization of
Algorithm \ref{alg:TwoLevelVCycle} is justified in practice but in
principle a detailed analysis of when Algorithm
\ref{alg:TwoLevelVCycle} becomes a direct solver is of clear interest
to the wider community. We also note that Algorithm
\ref{alg:TwoLevelVCycle} is more general than what is commonly used in
AMG literature in the sense of the smoothing parameters
$\alpha_i$. The standard choice is to set $\alpha_i = \alpha_0$ for
some $\alpha_0$, while in the following section we make specific and
varied choices of $\alpha_i$, $S^{-1}$ and $R$ to induce a direct
method.

Much of the analysis in this paper relies on error propagation
operators commonly used in multigrid analysis, see for example
\citep[eq. (2) and below]{MacLachlanManteuffelMcCormick2006}.  We
build the error propagation operator for Algorithm
\ref{alg:TwoLevelVCycle} piecemeal, starting with the error
propagation for a single smoothing step:
\begin{equation} \label{eq:Esgen}
  E^{(s)}(\alpha_i) = I - \alpha_i S^{-1} L.
\end{equation}
Likewise, the error propagation for the coarse space correction is
\begin{equation} \label{eq:Ecgen}
  E^{(c)} = I - P M_0^{-1} R L.
\end{equation}
The full error propagation operator for the algorithm is then
\begin{equation}\label{eqn:fullerror}
  E(\alpha_1, \dots, \alpha_m) =
  \left(\prod_{i=m}^{1} E^{(s)}(\alpha_i)\right) E^{(c)}
  \left(\prod_{i=1}^{m} E^{(s)}(\alpha_i)\right).
\end{equation}
It is linked with the preconditioned system through the equality
$M^{-1} L = I - E$.  Analysis of the eigenproperties of $E$ then
informs those of $M^{-1} L$.  In particular, clustering the spectrum
of $E$ around the origin will cluster that of $M^{-1}L$ around 1.

Combined with a nested dissection ordering and static condensation,
the algorithm is used by methods such as the Hierarchical
Poincaré--Steklov (HPS) method (see \citep{LuceroLorca2025,
  OutrataLuceroLorca2025} and references therein).

\section{New clustered method}\label{sec:newclustered}

\PL{
  We introduce a new algorithm, Algorithm
  \ref{alg:TwoLevelVCycleClustered}, based on Algorithm
  \ref{alg:TwoLevelVCycle} but with a nonsingular choice of the
  smoother
}
\begin{align}\label{eqn:optimalchoice2}
  \begin{aligned}
    \Pideal = \begin{bmatrix}
                -A^{-1} B \\
                I
              \end{bmatrix}, \quad
    \Rideal = \begin{bmatrix}
                -CA^{-1} & I
              \end{bmatrix}, \quad
    \Sideal = \begin{bmatrix}
                A^{-1} \\
                & D^{-1}
              \end{bmatrix} \\
    \implies \Mideal = \Rideal L \Pideal
    = D - C A^{-1} B,
  \end{aligned}
\end{align}
\PL{and a specific choice of the smoothing parameters}
\begin{equation}\label{eqn:optimalchoice3}
  \hat{\alpha}_i =
  \frac{1}{1 - \cos \left ( \frac{2 \pi i}{2m+1} \right )}.
\end{equation}

\begin{algorithm}
  \caption{TwoLevelVCycleClustered($\mathbf{v}$)}
  \begin{algorithmic}[1]
    \State $\mathbf{x} \leftarrow \mathbf{0}$
    \For{$i=1,\dots,m$}
    \State $\mathbf{x} \leftarrow \mathbf{x} +
    \hat{\alpha}_i \Sideal (\mathbf{v} - L \mathbf{x})$
    \EndFor
    \State $\mathbf{x} \leftarrow \mathbf{x} +
    \Pideal (\Mideal)^{-1} \Rideal (\mathbf{v} - L \mathbf{x})$
    \For{$i=m,\dots,1$}
    \State $\mathbf{x} \leftarrow \mathbf{x} +
    \hat{\alpha}_i \Sideal (\mathbf{v} - L \mathbf{x})$
    \EndFor
    \State \textbf{return} \PL{$\mathbf{x}$}
  \end{algorithmic}
  \label{alg:TwoLevelVCycleClustered}
\end{algorithm}

\PL{
  Going forward, we denote the matrix representation of
  Algorithm~\ref{alg:TwoLevelVCycleClustered} by $\hat{M}^{-1}$.  In
  order to analyze it, we study Algorithm~\ref{alg:TwoLevelVCycle}
  with the ideal operator choices~\eqref{eqn:optimalchoice2} but
  general smoothing parameters $\alpha_1,\dotsc ,\alpha_m$ and denote
  its matrix representation by $M^{-1}$.  With the choices of
  restriction, prolongation, and coarse space operators
  in~\eqref{eqn:optimalchoice2}, we evaluate the error propagation
  operators for the smoothing steps and coarse space correction from
  equations \eqref{eq:Esgen} and \eqref{eq:Ecgen} as
}
\begin{equation}\label{eq: Es}
  E^{(s)}(\alpha_i) = 
  \begin{bmatrix}
    (1-\alpha_i)I & -\alpha_i A^{-1}B\\
    -\alpha_i D^{-1}C& (1-\alpha_i)I
  \end{bmatrix},
\end{equation}
\begin{equation} \label{eq: Ec}
  E^{(c)} = 
  \begin{bmatrix}
    I & A^{-1}B\\
      & 
  \end{bmatrix}\PL{,}
\end{equation}
\PL{with the full error propagation operator given as
  in~\eqref{eqn:fullerror}.}

\PL{
  We note that $\Sideal$ acts as a block Jacobi smoother, using
  $A^{-1}$ and $D^{-1}$ instead of using only the former as
  in~\eqref{eqn:optimalchoice1}.  The ideal restriction and
  prolongation operators ensure that if $L$ is Hermitian, then so is
  $\Mideal$; in fact, the symmetric pre- and post-smoothing with
  $\Sideal$ extends this to $M^{-1}$, assuming we choose real
  smoothing coefficients, as we show next.
}

\PL{
  \begin{lemma}\label{lemma_secClstrMeth_Mhermitian}
    Let $\alpha_1,\dotsc ,\alpha_m \in \mathbb{R}$ and, in addition,
    assume that $L$ is Hermitian. Then the matrix representation of
    Algorithm~\ref{alg:TwoLevelVCycle} with the ideal choices
    in~\eqref{eqn:optimalchoice2} is Hermitian, i.e.,
    $M^{-1} = M^{-*}$.
  \end{lemma}
  \begin{proof}
    First, we transpose the claim to the language of the error
    operators~(\ref{eq: Es}--\ref{eq: Ec}): using the identity
    $M^{-1}L = I-E$, we have
    \begin{equation*}
      M^{-1} = M^{-*} 
      \quad \iff \quad
      L - L M^{-1} L = L - L M^{-*} L 
      \quad \iff \quad
      LE = E^{*}L,
    \end{equation*}
    with $E=E(\alpha_1, \dots, \alpha_m)$ given as
    in~\eqref{eqn:fullerror}. If $L$ is, moreover, HPD, then the
    above shows that $M^{-1}$ is Hermitian if and only if $E$ is
    $L$-Hermitian, i.e., $E$ is Hermitian with respect to the inner
    product induced by $L$. A direct calculation shows that
    \begin{equation*}
      \begin{aligned}
        L E^{(s)}(\alpha_i) = L - \alpha_i L\Sideal L = \left( E^{(s)}(\alpha_i) \right)^{*} L 
        \quad \text{ and } \quad 
        L E^{(c)} = L - L \Pideal \Mideal^{-1} \Rideal L = \left( E^{(c)} \right)^{*} L,
      \end{aligned}
    \end{equation*}
    so that, together with $L=L^{*}$ and~\eqref{eqn:fullerror}, we
    conclude that $LE = E^{*}L$, yielding the statement.
  \end{proof}
}

In Theorem \ref{theo:cluster} we prove that \PL{taking
  $\alpha_i=\hat{\alpha}_i$ as in~\eqref{eqn:optimalchoice3}}
\PL{clusters} the spectrum of $E(\alpha_1,\dots,\alpha_m)$ in two
points and \PL{makes} the operator diagonalizable, provided
$T=A^{-1}BD^{-1}C$ is diagonalizable. This then extends to the
preconditioned system, making it possible to obtain a direct method by
writing the minimal polynomial explicitly, or using Krylov subspace
methods. The next section provides details on the approach.

\section{New direct method}\label{sec:newdirect}

\PL{Provided} $T=A^{-1}BD^{-1}C$ is diagonalizable, we show later in
Theorem~\ref{theo:cluster} that \PL{$\hat{M}^{-1}L$} is diagonalizable
and has a two-point spectrum \PL{$\{1,1-\hat{\rho}_m\}$} with
\PL{$\hat{\rho}_m \equiv 1/(2m+1)^2$}. Therefore, the minimal polynomial of
\PL{$\hat{M}^{-1} L$} is
\begin{equation}\label{eqn_secDirMeth_MinvLMinimPoly}
  \left(\PL{\hat{M}^{-1}}L - I \right)
  \left(\PL{\hat{M}^{-1}}L - (\PL{1-}\hat{\rho}_m) I \right) = 0,
\end{equation}
and a rearrangement of the terms gives
\begin{equation}\label{eqn_secDirMeth_LinvAsFunctnOfLMinv}
  L^{-1} =
  \left( \left(1+\PL{\frac{1}{1-\hat{\rho}_m}}\right)
    I - \PL{\frac{1}{1-\hat{\rho}_m}} \PL{\hat{M}^{-1}}L \right) \PL{\hat{M}^{-1}}.
\end{equation}
We have successfully constructed a direct two-level method; it
requires two applications of
Algorithm~\ref{alg:TwoLevelVCycleClustered} and a single application
of $L$ in the matrix-vector product sense. In this sense, our choices
of $S^{-1}, P, R, M_0$ and
\PL{$\hat{\alpha}_1,\dotsc , \hat{\alpha}_m$} are optimal, analogous
to the optimality of the choices in eq. \eqref{eqn:optimalchoice1} for
Algorithm~\ref{alg:TwoLevelVCycle}.  \PL{First, we highlight that when
  applying $\hat{M}^{-1}$}
\PL{in~}\eqref{eqn_secDirMeth_LinvAsFunctnOfLMinv}, we are free to
retrieve \PL{$\Mideal^{-1}$} in
Algorithm~\ref{alg:TwoLevelVCycleClustered} by using \PL{the
  formula~}\eqref{eqn_secDirMeth_LinvAsFunctnOfLMinv} again. The
method can be applied recursively as desired.  \PL{Second, a natural
  step is to use $\hat{M}^{-1}$ as a preconditioner; we start by
  showing that when $L$ is HPD, so is $\hat{M}^{-1}$. To do that, we
  are going to use the result of Theorem~\ref{theo:cluster} that
  $\hat{M}^{-1}L$ is diagonalizable and has a two-point spectrum
  $\{1,1-\hat{\rho}_m\}$; Proposition~\ref{prop_secDirrMeth_Mhpd} is
  not used to obtain results of Theorem~\ref{theo:cluster} and we
  state it here, before Theorem~\ref{theo:cluster} only for the sake
  of exposition.}

\PL{
\begin{proposition}\label{prop_secDirrMeth_Mhpd}
  If $L$ is HPD and $T:=A^{-1}BD^{-1}C$ is diagonalizable, then the
  matrix representation of
  Algorithm~\ref{alg:TwoLevelVCycleClustered}, i.e., $\hat{M}^{-1}$,
  is also HPD.
\end{proposition}
\begin{proof}
  Thanks to Lemma~\ref{lemma_secClstrMeth_Mhermitian} it remains to
  show the positive definiteness of $\hat{M}^{-1}$.
  Theorem~\ref{theo:cluster} gives the spectrum of $\hat{M}^{-1}L$ as
  the two points $\{1,1-\hat{\rho}_{m}\}$, both positive real
  numbers. Moreover, the similarity transformation
  \begin{equation*}
    \hat{M}^{-1}L \sim  L^{1/2} \hat{M}^{-1} L^{1/2}
  \end{equation*}
  shows that the spectra of these matrices coincide. The statement
  follows by noticing that $L^{1/2} \hat{M}^{-1} L^{1/2}$ is
  positive-definite if and only if $\hat{M}^{-1}$ is.
\end{proof}
}

\PL{Proposition~\ref{prop_secDirrMeth_Mhpd} enables the use of
  $\hat{M}^{-1}$ as a preconditioner for CG when $L$ is HPD; notice
  that $\hat{M}^{-1}L$ itself need not be HPD; it is enough that both
  $L$ and $\hat{M}^{-1}$ are. This is because preconditioned CG is
  equivalent to the standard CG running with respect to the inner
  product induced by the preconditioner, in this case $\hat{M}^{-1}$,
  see~\citep[Section 9.2]{saad2003iterative} and
  also~\citep{malek2014preconditioning}.}  Using Algorithm
\ref{alg:TwoLevelVCycleClustered} to precondition CG (in the HPD case)
or GMRES we have the following proposition.
\begin{proposition}\label{ref:KrylovConvergence}
  If the matrix $E$ in~\eqref{eqn:fullerror} for Algorithm
  \ref{alg:TwoLevelVCycleClustered} is diagonalizable and its spectrum
  is clustered into two points, then both CG (assuming $L$ is HPD) and
  GMRES converge in \PL{at most two iterations}.
\end{proposition}
\begin{proof}
  The proof is obtained simply by realizing that the minimal
  polynomial of a diagonalizable matrix with a spectrum clustered in
  two eigenvalues is of degree 2, and thus both CG and GMRES converge
  in at most two iterations, see \citep[Theorems 2.2.3 and
  2.3.1]{LiesenStrakos2013}.
\end{proof}

\begin{remark}\label{cor:KSMinstdCGGMRES}
  The results of Proposition \ref{ref:KrylovConvergence} hold for
  \PL{other} optimal Krylov subspace methods such as
  MINRES~\citep[Section 2.3]{LiesenStrakos2013}, not only CG or
  GMRES. However, they cannot be directly extended to non-optimal
  Krylov subspace methods, such as FOM, QMR or BiCG,
  see~\citep[Section 2.5.6]{LiesenStrakos2013}. In what follows, we
  abbreviate the optimal Krylov subspace methods by
  $\mathrm{KSM}^{\star}$.
\end{remark}

With Proposition \ref{ref:KrylovConvergence}, we can use Algorithm
\ref{alg:TwoLevelVCycleClustered} to precondition a
$\mathrm{KSM}^{\star}$ method and recursively solve the coarse problem
using the same algorithm (see \citep{NotayVassilevski2008}), this is
referred to as the K-cycle. We write
$\mathrm{KSM}^{\star}(L,\mathbf{b},P^{-1})$ to refer to using an
optimal Krylov subspace method (such as CG or GMRES) to solve the
linear system $P^{-1}L\mathbf{x}=P^{-1}\mathbf{b}$ with $P^{-1}$ as
the left preconditioner, and define the direct AMG algorithm below.

\begin{algorithm}
  \caption{AMGK($\mathbf{v}$)}
  \begin{algorithmic}[1]
    \State $\mathbf{x} \leftarrow \mathbf{0}$
    \For{$i=1,\dots,m$}
    \State $\mathbf{x} \leftarrow \mathbf{x} +
    \hat{\alpha}_i \Sideal (\mathbf{v} - L \mathbf{x})$
    \EndFor
    \If {we have reached the coarsest level}
    \State $\mathbf{x} \leftarrow \mathbf{x} +
    \Pideal \Mideal^{-1} \Rideal (\mathbf{v} - L \mathbf{x})$
    \Else
    \State $\mathbf{x} \leftarrow \mathbf{x} +
    \Pideal
    \bigg(\mathrm{KSM^{\star}}\PL{\bigg(
    \Mideal,
    \Rideal (\mathbf{v} - L \mathbf{x}),
    \mathrm{AMGK}(\cdot)
    \bigg)}\bigg)$
    \EndIf
    \For{$i=m,\dots,1$}
    \State $\mathbf{x} \leftarrow \mathbf{x} +
    \hat{\alpha}_i \Sideal (\mathbf{v} - L \mathbf{x})$
    \EndFor
    \State \textbf{return} \PL{$\mathbf{x}$}
  \end{algorithmic}
  \label{alg:AMGK}
\end{algorithm}

\begin{lemma}\label{lem:AMGK}
  If Algorithm \ref{alg:AMGK} is used to precondition
  $\mathrm{KSM}^{\star}$, every $\mathrm{KSM}^{\star}$ call in the
  hierarchy converges to the exact solution in 2 iterations in exact
  arithmetic.
\end{lemma}
\begin{proof}
  Let $\ell_0$ be the coarsest level in the recursion, where we choose
  to invert $\Mideal$ directly.

  Consider the immediate finer level $\ell_1$. Since we have chosen
  the operators and constants from eq.  \eqref{eqn:optimalchoice2},
  the spectrum of the preconditioned operator is clustered in two
  values and the minimal polynomial is of degree 2, thus the
  $\mathrm{KSM}^{\star}$ converges in two iterations in exact
  arithmetic.

  Levels $\ell_1$ and $\ell_0$ are now a direct method from the point
  of view of $\ell_2$, since the $\mathrm{KSM}^{\star}$ converges in
  two iterations. The result follows by induction on $\ell_i$, for
  $i=2$ until the finest level.
\end{proof}

\begin{remark}
  Analogous results hold also when we consider right-preconditioning
  as the default.
\end{remark}

The next sections give a complete description of the spectrum of
Algorithm \ref{alg:TwoLevelVCycleClustered} and a proof of the
spectral clustering claimed in this section.

\section{Spectral analysis of $E$} \label{sec:spectralanalysis}

This section analyzes the spectral properties of the error propagation
operator $E(\alpha_1,\dotsc ,\alpha_m)$ of \PL{Algorithm
\ref{alg:TwoLevelVCycle} with the ideal operator
choices~\eqref{eqn:optimalchoice2} for general smoothing coefficients
$\alpha_1,\dotsc ,\alpha_m$.}  Related ideas can be found in~\citep[\S
2.4]{Hackbusch1985} for geometric multigrid, and in
\citep{SiefertdeSturler2006,Ipsen2001} for saddle-point systems.

Consider the characteristic polynomial
$\det \left (t I - E^{(s)}(\alpha_i) \right )$ of $E^{(s)}(\alpha_i)$
for some $\alpha_i \neq 0$. It is straightforward to show that
\begin{equation}\label{eqn_secSpctrlAnal_EsEigvalsParametrztn}
  \det \left (t I - E^{(s)}(\alpha_i) \right ) =
  \PL{(\alpha_i\mu)^{2n-N}}
  \det\left( (\PL{\alpha_i} \mu)^2 I_{N-n} - \PL{\alpha_i^2} T \right),
\end{equation}
where $T := A^{-1}B D^{-1}C$ and $\mu := (t + \alpha_i - 1)/\alpha_i$.
Thus, the eigenvalues of the error propagation operator for smoothing
are parametrized by those of $T$ and by $\alpha_i$. In this section we
show that in fact the eigenpairs of the error propagation operator
$E(\alpha_1,\dotsc ,\alpha_m)$ itself can be explicitly parametrized
by the eigenvalues of $T$ and by $\alpha_1, \dotsc , \alpha_m$.

For convenience, we will denote the nonzero eigenvalues of $T$ by
$\lambda^2$, where $\lambda$ coincides with $\mu$ at the roots of the
characteristic polynomial in
eq. \eqref{eqn_secSpctrlAnal_EsEigvalsParametrztn}\footnote{So that
  the symbol $\lambda$ corresponds to the principal branch of the
  square root, i.e.,
  $\lambda = \exp \left(\tfrac{1}{2}\log \lambda^2 \right)$, where
  $\log$ denotes the principal complex logarithm, with branch cut on
  $(-\infty,0]$ and \PL{$\arg \lambda \in (-\pi/2,\pi/2]$}.}. This
does not reflect any assumption on symmetry or definiteness of $T$;
rather, similarly to~\eqref{eqn_secSpctrlAnal_EsEigvalsParametrztn},
the computation naturally induces this notation. Equivalently, we
could label the eigenvalues of $T$ by $\lambda$ and replace $\lambda$
by $\sqrt{\lambda}$ in the analysis below. We note that the entire
analysis is independent of the sign choice for $\lambda$ as the
resulting expressions are either even in $\lambda$ or involve only the
symmetric combination of $\lambda$ and $-\lambda$.

With this notation, we would also like to highlight that in some
applications it is more natural to think of an eigenpair
$(\lambda^2,\mathbf{w})$ of $T$ as a generalized eigenpair of the
Schur complement matrix pencil $\{B D^{-1} C, A\}$,
\begin{equation*}
(B D^{-1} C) \mathbf{w} = \lambda^2 A \mathbf{w},
\end{equation*}
as the generalized eigenpairs can have a concrete interpretation in
the original problem. Here we will work directly with $T$.

\begin{remark}\label{rmk_secSpctrlAnal_Notation}
  In the work that follows, we make two assumptions:
  \begin{itemize}
  \item $n \geq N/2$;
  \item the matrix $T$ is diagonalizable.
  \end{itemize}
  Through notation we distinguish two types of eigenpairs of the matrix $T$:
  \begin{itemize}
  \item eigenpairs with nonzero eigenvalues, denoted by $(\lambda^2,\mathbf{w})$;
  \item eigenpairs with zero eigenvalues, denoted by $(0,\mathbf{z})$,
  \end{itemize}
  so that in our notation  $\lambda^2\neq 0$.
\end{remark}

We note that the case $n<N/2$ can be handled by symmetry.  In the
calculation above, we could use the other Schur complement formulation
for the determinant, featuring the matrix $D^{-1}CA^{-1}B$ instead of
$T$.  This would lead to calculations that are analogous to the ones
we carry out below, replacing the eigenpairs of $T$ with those of
$D^{-1}CA^{-1}B$ and exchanging the blocks $D^{-1}C$ and $A^{-1}B$
wherever necessary.

For each type of eigenpair of $T$, we define a matrix whose column
space is a shared invariant subspace of $E^{(c)}$ and
$E^{(s)}(\alpha_i)$ for all $i$, in turn giving invariant subspaces of
$E$. Through these we identify eigenpairs of $E$ corresponding to
nonzero eigenvalues of $E$. \PL{We then show that
  $\mathbb{C}^N = \mathcal{V} \oplus \mathcal{W}$, where $\mathcal{V}$
  is the direct sum of these invariant subspaces and
  $\mathcal{W} \subset \Ker$
  (Proposition~\ref{prop_SpctrlAnal_InvSbspc}), and, in
  Section~\ref{sec:spectralclustering}, that for
  $\alpha_i = \hat{\alpha}_i$ the operator $E$ is diagonalizable
  (Theorem~\ref{theo:cluster}).}

The calculations share some similarities to local Fourier analysis
(LFA; see~\citep[Section 4]{LuceroLorcaGander2024}), but our
derivations are fundamentally different. In contrast to LFA we do not
assume periodicity, i.e., we work with a general matrix rather than
with a circulant one. While $L$ may in fact correspond to a
discretization of some infinite-dimensional operator, the analysis
here applies no matter the choice of the operator or the
discretization, provided $T$ is diagonalizable.

\subsection{Nonzero eigenvalues}

For each nonzero eigenvalue $\lambda^2$, define the matrix
$V_\lambda \in \mathbb{C}^{N \times 2}$ as
\begin{equation} \label{eq: non zero subspace}
  V_\lambda := \begin{bmatrix}
                 \lambda \vec{w} & \lambda \vec{w} \\
                 D^{-1} C \vec{w} & -D^{-1} C \vec{w}
               \end{bmatrix},
\end{equation}
and the two-dimensional subspace
$\mathcal{V}_\lambda \subset \mathbb{C}^N$ as the column space of
$V_\lambda$.  The action of $E^{(s)}(\alpha_i)$ on $V_\lambda$ is
\begin{equation*}
  \begin{aligned}
    E^{(s)}(\alpha_i) V_{\lambda} &= 
    \begin{bmatrix}
      \left( (1-\alpha_i) \lambda - \alpha_i \lambda^2 \right)\mathbf{w} & \PL{\left((1-\alpha_i)\lambda + \alpha_i\lambda^2 \right)} \mathbf{w} \\
      \left( 1- \alpha_i -\alpha_i\lambda \right)D^{-1}C\mathbf{w} & -\left( 1- \alpha_i + \alpha_i\lambda \right)D^{-1}C\mathbf{w}
    \end{bmatrix} \\
    &= V_\lambda \begin{bmatrix}
                   1 - \alpha_i - \alpha_i \lambda \\
                   & 1 - \alpha_i + \alpha_i \lambda
                 \end{bmatrix},
  \end{aligned} 
\end{equation*}
and the action of $E^{(c)}$ is
\begin{equation*}
  E^{(c)} V_{\lambda} = 
  \begin{bmatrix}
    I & A^{-1}B \\
      &
  \end{bmatrix}
  V_{\lambda} = 
  \frac{1}{2} V_\lambda
  \begin{bmatrix}
    1 + \lambda & 1 - \lambda \\
    1 + \lambda & 1 - \lambda
  \end{bmatrix} =
  \frac{1}{2} V_\lambda
  \begin{bmatrix}
    1 \\
    1
  \end{bmatrix}
  \begin{bmatrix}
    1 + \lambda & 1 - \lambda
  \end{bmatrix}.
\end{equation*}
In words, for every nonzero eigenvalue of $T$ we have found a
two-dimensional subspace $\mathcal{V}_{\lambda} \subset \mathbb{C}^N$
that is invariant under both $E^{(c)}$ and $E^{(s)}(\alpha_i)$ for any
$\alpha_i \in \mathbb{C}$. Since the full error propagation operator
is a composition of these, we have that $\mathcal{V}_{\lambda}$ is
also an invariant subspace of $E(\alpha_1, \dots,
\alpha_m)$. Moreover, we have represented each of the operators
$E^{(c)}, E^{(s)}(\alpha_i)$ and $E(\alpha_1, \dots, \alpha_m)$
restricted to $\mathcal{V}_{\lambda}$ as a $2 \times 2$ matrix; we
denote these by $E_{\lambda}^{(c)}, E_{\lambda}^{(s)}(\alpha_i)$ and
$E_{\lambda}(\alpha_1, \dots, \alpha_m)$ with
\begin{equation*}
  E_{\lambda}^{(c)} := \frac{1}{2}
  \begin{bmatrix}
    1 \\
    1
  \end{bmatrix}
  \begin{bmatrix}
    1 + \lambda & 1 - \lambda
  \end{bmatrix}, \quad
  E_{\lambda}^{(s)}(\alpha_i) :=
  \begin{bmatrix}
    1 - \alpha_i - \alpha_i \lambda \\
    & 1 - \alpha_i + \alpha_i \lambda
  \end{bmatrix}
\end{equation*} 
and $E_{\lambda}(\alpha_1, \dots, \alpha_m)$ defined analogously
to~\eqref{eqn:fullerror} as
\begin{equation*}
  E_{\lambda}(\alpha_1, \dots, \alpha_m) :=
  \left(\prod_{i=m}^{1} E_{\lambda}^{(s)}(\alpha_i)\right)
  E_{\lambda}^{(c)}
  \left(\prod_{i=1}^{m} E_{\lambda}^{(s)}(\alpha_i)\right) \equiv
  E_{\lambda}^{(s)} E_{\lambda}^{(c)} E_{\lambda}^{(s)},
\end{equation*}
where we introduce the notation for the product of the diagonal
matrices\footnote{Note that since the restriction of
  $E^{(s)}(\alpha_i)$ to $\mathcal{V}_\lambda$ is diagonal, the order
  of the products is arbitrary.} $E_{\lambda}^{(s)}(\alpha_i)$
\begin{equation*}
  E^{(s)}_\lambda := \prod_{i=1}^m
  \begin{bmatrix}
    1 - \alpha_i - \alpha_i \lambda \\
    & 1 - \alpha_i + \alpha_i \lambda
  \end{bmatrix} =
  \begin{bmatrix}
    \rho_- \\
    & \rho_+
  \end{bmatrix},
  \quad \rho_\pm =
  \prod_{i=1}^m (1 - \alpha_i \pm \alpha_i \lambda).
\end{equation*}
By definition, eigenpairs of the operators
$E^{(c)}, E^{(s)}(\alpha_i)$ and, most importantly,
$E(\alpha_1, \dots, \alpha_m)$ can be retrieved from those of
$E_{\lambda}^{(c)}, E_{\lambda}^{(s)}(\alpha_i)$ and
$E_{\lambda}(\alpha_1, \dots, \alpha_m)$. \PL{ For the moment, assume
  that $\rho_{\pm}\neq 0$ and $E_{\lambda}^{(s)}$ is nonsingular. Then
  \begin{align}\label{eq: calc nnz eigvec LambdaNnz}
    E_\lambda \frac{1}{2} E^{(s)}_\lambda
    \begin{bmatrix}
      1 \\
      1
    \end{bmatrix}
    & = \frac{1}{2} E^{(s)}_\lambda
      \begin{bmatrix}
        1 \\
        1
      \end{bmatrix}
      \left (
      \begin{bmatrix}
        1 + \lambda & 1 - \lambda
      \end{bmatrix}
      E^{(s)}_\lambda \frac{1}{2} E^{(s)}_\lambda
      \begin{bmatrix}
        1 \\
        1
      \end{bmatrix}
      \right ),
  \end{align}
i.e.,} the vector
$(1/2) E^{(s)}_\lambda \begin{bmatrix} 1 & 1 \end{bmatrix}^T$ is an
eigenvector of $E_\lambda$ with the corresponding eigenvalue
\begin{equation}\label{eq:rho}
\begin{aligned}
  \rho_{\PL{\lambda^2}} :=& \frac{1}{2}
  \begin{bmatrix}
    1 + \lambda & 1 - \lambda
  \end{bmatrix}
  \left(
  E^{(s)}_\lambda \right )^2
  \begin{bmatrix}
    1 \\
    1
  \end{bmatrix} =
  \frac{1}{2} \left(
  (1 + \lambda) \rho_-^2 + (1 - \lambda) \rho_+^2 \right ) \\
  =& \frac{1}{2} (1 + \lambda)
     \prod_{i=1}^m (1 - \alpha_i - \alpha_i \lambda)^2 +
     \frac{1}{2} (1 - \lambda)
     \prod_{i=1}^m (1 - \alpha_i + \alpha_i \lambda)^2.
\end{aligned}
\end{equation}
By the same argument, the second eigenvalue of $E_\lambda$ is zero
with eigenvector $\vec{y}$ such that
\begin{align*}
  \begin{bmatrix}
    1 + \lambda & 1 - \lambda
  \end{bmatrix}
  E^{(s)}_\lambda \vec{y} = 0
  \implies \vec{y} =
  \begin{bmatrix}
    \rho_+ (1 - \lambda) \\
    -\rho_- (1 + \lambda)
  \end{bmatrix}.
\end{align*}

\PL{
  \begin{remark}\label{rmk_SpctrlAnal_rhoLambdaSqrd}
    Notice that $\rho_{\PL{\lambda^2}}$ is a polynomial in $\lambda$ that
    is symmetrical, i.e., takes the same value at $\lambda$ and
    $-\lambda$. Hence, it is in fact a polynomial in $\lambda^2$,
    justifying the subscript notation.
  \end{remark}
}
\PL{
  If $E_{\lambda}^{(s)}$ becomes singular, the computations can be
  directly adjusted; we summarize this in the following proposition.
}

\PL{
  \begin{proposition} \label{prop_SpctrlAnal_NonZeroEigvlsT} 
    Let $(\lambda^2, \vec{w})$ be an eigenpair of $T$ with
    $\lambda^2 \neq 0$ and fix
    $\alpha_1,\dotsc ,\alpha_m \in \mathbb{C}$. If $\rho_{\pm}\neq 0$,
    then $E(\alpha_1, \dots, \alpha_m)$ has the two eigenpairs
    \begin{equation*}
      \left ( 0,
        \begin{bmatrix}
          \left ( \rho_+ (1-\lambda) - \rho_- (1+\lambda) \right ) \lambda \vec{w} \\
          \left ( \rho_+ (1-\lambda) + \rho_- (1+\lambda) \right ) D^{-1} C \vec{w}
        \end{bmatrix} \right ),
      \quad
      \left ( \rho_{\PL{\lambda^2}}, \frac{1}{2}
        \begin{bmatrix}
          (\rho_- + \rho_+) \lambda \vec{w} \\
          (\rho_- - \rho_+) D^{-1} C \vec{w}
        \end{bmatrix} \right ).
    \end{equation*}
    If $\rho_{+} = 0$ while $\rho_{-} \neq 0$, then
    $E(\alpha_1, \dots, \alpha_m)$ has the two eigenpairs
    \begin{equation*}
      \left( 0,
        \begin{bmatrix}
          \lambda \vec{w} \\
          -D^{-1}C\vec{w}
        \end{bmatrix}
      \right), \quad
      \left( \rho_{\PL{\lambda^2}},
        \begin{bmatrix}
          \lambda \vec{w} \\
          D^{-1}C\vec{w}
        \end{bmatrix} \right).
    \end{equation*}
    If $\rho_{-}= 0$ while $\rho_{+} \neq 0$, then
    $E(\alpha_1, \dots, \alpha_m)$ has the two eigenpairs
    \begin{equation*}
      \left( 0,
        \begin{bmatrix}
          \lambda \vec{w} \\
          D^{-1}C\vec{w}
        \end{bmatrix} \right), \quad
      \left( \rho_{\PL{\lambda^2}},
        \begin{bmatrix} 
          \lambda \vec{w}\\
          -D^{-1}C\vec{w} 
        \end{bmatrix} \right).
    \end{equation*}
    If both $\rho_{+} = \rho_{-} = 0$, \PL{then 
      $E(\alpha_1, \dots, \alpha_m) V_\lambda = 0$; in particular the columns of $V_{\lambda}$ are eigenvectors 
      of $E(\alpha_1, \dots, \alpha_m)$ corresponding to the eigenvalue $0$.}
  \end{proposition}
  \begin{proof}
    If $\rho_{\pm}\neq 0$, then the calculations in~\eqref{eq: calc
      nnz eigvec LambdaNnz} and below give the eigenpairs of
    $E_{\lambda}$ and writing with respect to the basis~\eqref{eq: non
      zero subspace} gives the first part of the statement.  Assuming
    $\rho_{+}=0$, we observe that
    \begin{equation*}
      \frac{1}{2} E_{\lambda}^{(s)}
      \begin{bmatrix}
        1 \\
        1
      \end{bmatrix} =
      \frac{\rho_{-}}{2}
      \begin{bmatrix}
        1 \\
        0
      \end{bmatrix}
      \quad \mathrm{and} \quad
      E_{\lambda}^{(s)}
      \begin{bmatrix}
        0 \\
        1
      \end{bmatrix} =
      \begin{bmatrix}
        0 \\
        0
      \end{bmatrix}.
    \end{equation*}
    First, as~\eqref{eq: calc nnz eigvec LambdaNnz} does not require
    $\rho_{+}\neq 0$, we conclude that $\begin{bmatrix} 1 & 0 \end{bmatrix}^T$ is an
    eigenvector of $E_{\lambda}$ corresponding to the eigenvalue
    $\rho_{\PL{\lambda^2}} = (1+\lambda)/2 \cdot \rho_{-}^2$. Next, the
    second equality shows that
    $\begin{bmatrix} 0 & 1 \end{bmatrix}^T$ lies in the
    kernel of $E^{(s)}_{\lambda}$ and hence of $E_{\lambda}$. Writing
    the vectors $\begin{bmatrix} 1 & 0 \end{bmatrix}^T$,
    $\begin{bmatrix} 0 & 1 \end{bmatrix}^T$ with respect to the basis~\eqref{eq: non
      zero subspace} yields the result. The case $\rho_{-} = 0$ is
    fully analogous \PL{to that of $\rho_{+}=0$.  Having both}
    $\rho_{\pm} = 0$ implies 
    \PL{$E(\alpha_1, \dots, \alpha_m) V_\lambda = 0$ and hence the statement.}
  \end{proof}
}

Next, we adapt these ideas also for the case of eigenvectors
corresponding to zero eigenvalues of $T$.

\subsection{Zero eigenvalues}
We start by taking an eigenvector $\vec{z}$ of $T$ with eigenvalue
zero such that $\vec{z} \notin\mathrm{Ker}(C)$ and address the case
$\vec{z} \in\mathrm{Ker}(C)$ later in this section. We set
$V_{\mathbf{z}} \in \mathbb{C}^{N \times 2}$ as
\begin{equation}\label{eq: Vbasis LambdaZero}
  V_{\mathbf{z}} :=
  \begin{bmatrix}
    & \vec{z} \\
    D^{-1} C \vec{z}
  \end{bmatrix},
\end{equation}
and the subspace $\mathcal{V}_{\mathbf{z}} \subset \mathbb{C}^N$ as
the column space of $V_{\mathbf{z}}$.  The action of
$E^{(s)}(\alpha_i)$ and $E^{(c)}$ on $V_{\mathbf{z}}$ then becomes
\begin{equation*}
  E^{(s)}(\alpha_i) V_{\mathbf{z}} = V_{\mathbf{z}}
  \begin{bmatrix}
    1 - \alpha_i & -\alpha_i \\
                 & 1 - \alpha_i
  \end{bmatrix} = V_{\mathbf{z}} \left ( I - \alpha_i
    \begin{bmatrix}
      1 & 1 \\ & 1
    \end{bmatrix}
  \right ), \quad
  E^{(c)} V_{\mathbf{z}} = V_{\mathbf{z}}
  \begin{bmatrix} ~ \\
    & 1
  \end{bmatrix},
\end{equation*}
i.e., $\mathcal{V}_{\mathbf{z}}$ is an invariant subspace of
$E^{(s)}(\alpha_i)$ and $E^{(c)}$ and thus also of
$E(\alpha_1, \dots, \alpha_m)$. Analogously to the previous section,
we set
\begin{equation*}
  E_{\mathbf{z}}^{(c)} :=
  \begin{bmatrix} ~ \\
    &
      1
  \end{bmatrix}, \quad
  E_{\mathbf{z}}^{(s)}(\alpha_i) := I -
  \alpha_i
  \begin{bmatrix}
    1 & 1 \\
      & 1
  \end{bmatrix}, \quad
  E_{\mathbf{z}}^{(s)} := \prod_{i=1}^m \left ( I -
    \alpha_i
    \begin{bmatrix}
      1 & 1 \\
        & 1
    \end{bmatrix} \right ).
\end{equation*} 
Although $E_{\mathbf{z}}^{(s)}(\alpha_i)$ are no longer diagonal, we
see that these still commute for any $\alpha_i$ and so we can again
obtain eigenpairs of $E(\alpha_1, \dots, \alpha_m)$ via calculating
those of $E_{\mathbf{z}}(\alpha_1, \dots, \alpha_m)$ with
\begin{equation*}
  E_{\mathbf{z}}(\alpha_1, \dots, \alpha_m) :=
  E_{\mathbf{z}}^{(s)} E_{\mathbf{z}}^{(c)} E_{\mathbf{z}}^{(s)}.
\end{equation*}
To that end, we introduce the elementary symmetric polynomials,
\begin{equation*}
  e_1(\alpha_1, \dots, \alpha_m) = \sum_{i=1}^m \alpha_i, \quad
  e_k(\alpha_1, \dots, \alpha_m) = \sum_{i=1}^{m-k+1}
  \alpha_i e_{k-1}(\alpha_{i+1}, \dots, \alpha_m),
\end{equation*}
which give us a convenient formula for $E^{(s)}_{\mathbf{z}}$:
\begin{equation*}
  E^{(s)}_{\mathbf{z}} = \prod_{i=1}^m \left ( I - \alpha_i
    \begin{bmatrix}
      1 & 1 \\
        & 1
    \end{bmatrix} \right ) =
  I + \sum_{k=1}^m (-1)^k e_k(\alpha_1, \dots, \alpha_m)
  \begin{bmatrix}
    1 & k \\
      & 1
  \end{bmatrix}.
\end{equation*}
A direct calculation then gives
\begin{align}\label{eq: calc eigvec LambdaZero}
  \begin{aligned}
    E_{\mathbf{z}}^{(c)} E^{(s)}_{\mathbf{z}} & =
    \begin{bmatrix}
      ~ \\
      & 1
    \end{bmatrix}
    \left ( 1 + \sum_{k=1}^m (-1)^k e_k(\alpha_1, \dots, \alpha_m) \right ) =
    E_{\mathbf{z}}^{(c)} \prod_{i=1}^m (1 - \alpha_i), \\
    E_{\mathbf{z}}^{(s)} E_{\mathbf{z}}^{(c)} & =
    \begin{bmatrix}
      0 & \sum\limits_{k=1}^m (-1)^k k e_k(\alpha_1, \dots, \alpha_m) \\
      0 & 1 + \sum\limits_{k=1}^m (-1)^k e_k(\alpha_1, \dots, \alpha_m)
    \end{bmatrix} =
    \begin{bmatrix}
      0 & \sum\limits_{i=1}^m (-\alpha_i) \prod\limits_{j \neq i} (1 - \alpha_j) \\
      0 & \prod\limits_{i=1}^m (1 - \alpha_i)
    \end{bmatrix} , \\
    \implies E_{\mathbf{z}} & =
    \begin{bmatrix}
      0 & \sum\limits_{i=1}^m (-\alpha_i) (1 - \alpha_i) \prod\limits_{j \neq i} (1 - \alpha_j)^2 \\
      0 & \prod\limits_{i=1}^m (1 - \alpha_i)^2
    \end{bmatrix}.
  \end{aligned}
\end{align}
The eigenvectors of $E_{\mathbf{z}}$ are then
$\begin{bmatrix} 1 & 0 \end{bmatrix}^T$ and the rightmost column of
$E_{\mathbf{z}}$, with corresponding eigenvalues 0 and
$\prod_{i=1}^m (1 - \alpha_i)^2$.

\begin{proposition} \label{prop_SpctrlAnal_ZeroEigvlsT}
  \PL{Let $(0, \vec{z})$ be} an eigenpair of $T$ such that
  $\vec{z} \notin \mathrm{Ker}(C)$ \PL{and fix
    $\alpha_1,\dotsc ,\alpha_m \in \mathbb{C}$.  If
    $\alpha_1, \dots, \alpha_m \neq 1$, then
    $E(\alpha_1, \dots, \alpha_m)$ has the two eigenpairs}
  \begin{equation*}
    \left ( 0,
      \begin{bmatrix}
        \vec{0} \\ D^{-1} C \vec{z}
      \end{bmatrix} \right )
    \quad \& \quad
    \left ( \prod_{i=1}^m (1 - \alpha_i)^2,
      \begin{bmatrix}
        \prod_{i=1}^m (1-\alpha_i) \vec{z} \\
        \sum_{i=1}^m (-\alpha_i) \prod_{j \neq i} (1 - \alpha_j) D^{-1} C \vec{z}
      \end{bmatrix} \right ).
  \end{equation*}
  \PL{Otherwise, $E(\alpha_1, \dots, \alpha_m)$ has the two
    eigenpairs}
  \begin{equation*}
    \left( 0, \left[
        \begin{matrix}
          \vec{0} \\ D^{-1}C\vec{z}
        \end{matrix}
      \right] \right) \quad \& \quad
    \left( 0, \left[
        \begin{matrix} 
          \vec{z} \\
          \vec{0}
        \end{matrix}\right] \right).
  \end{equation*}
\end{proposition}
\PL{
  \begin{proof}
    Notice the calculations leading to and including~\eqref{eq: calc
      eigvec LambdaZero} hold regardless of
    $\alpha_1,\dotsc ,\alpha_m$. The first part of the statement
    follows from writing the calculated eigenvectors below~\eqref{eq:
      calc eigvec LambdaZero} with respect to the basis~\eqref{eq:
      Vbasis LambdaZero}. However, if
    $\prod_{i=1}^m (1-\alpha_i) = 0$, then the calculated eigenvectors
    become identical (up to scaling); returning to~\eqref{eq: calc
      eigvec LambdaZero} we observe that $E_{\vec{z}}=0$ and hence has
    the eigenpairs $(0,\begin{bmatrix} 1 & 0 \end{bmatrix}^T)$ and
    $(0,\begin{bmatrix} 0 & 1 \end{bmatrix}^T)$ and writing these with
    respect to the basis~\eqref{eq: Vbasis LambdaZero} finalizes the
    proof.
  \end{proof}
}

Special care needs to be taken for eigenvectors $\mathbf{z}$ of $T$
that lie in the kernel of $C$. We see that the \PL{first} column of
$V_{\mathbf{z}}$ zeros out and, more importantly, a direct calculation
shows that its \PL{second} column becomes an eigenvector of
$E_{\mathbf{z}}$.

\begin{proposition} \label{prop_SpctrlAnal_ZeroEigvlsT_KerC} If
  $(0,\vec{z})$ is an eigenpair of $T$ such that
  $\vec{z} \in \mathrm{Ker}(C)$, then
  \begin{equation*}
    \left ( \prod_{i=1}^m (1 - \alpha_i)^2 ,
      \begin{bmatrix}
        \vec{z} \\
        ~
      \end{bmatrix} \right )
  \end{equation*}
  is an eigenpair of $E(\alpha_1, \dots, \alpha_m)$.
\end{proposition}
\begin{proof}
\begin{equation*}
  E^{(s)}(\alpha_i)
  \begin{bmatrix}
    \vec{z} \\
    ~
  \end{bmatrix} = (1-\alpha_i)
  \begin{bmatrix}
    \vec{z} \\
    ~
  \end{bmatrix},
  \quad E^{(c)}
  \begin{bmatrix}
    \vec{z} \\
    ~
  \end{bmatrix} =
  \begin{bmatrix}
    \vec{z} \\
    ~
  \end{bmatrix},
  \implies
  E
  \begin{bmatrix}
    \vec{z} \\
    ~
  \end{bmatrix} = \prod_{i=1}^m (1 - \alpha_i)^2
  \begin{bmatrix}
    \vec{z} \\
    ~
  \end{bmatrix}.
\end{equation*}
\end{proof}

This shows that the subspace $\mathcal{V}_{\mathbf{z}}$ formed as the
column space of $V_{\mathbf{z}}$ is still an invariant subspace of
$E$, but only one-dimensional.

\begin{remark}\label{rmk: rho_lambda}
  Notice that the nonzero eigenvalues of
  $E(\alpha_1, \dots, \alpha_m)$ from
  Propositions~\ref{prop_SpctrlAnal_ZeroEigvlsT}
  and~\ref{prop_SpctrlAnal_ZeroEigvlsT_KerC} correspond to the nonzero
  eigenvalue $\rho_{\PL{\lambda^2}}$ from
  Proposition~\ref{prop_SpctrlAnal_NonZeroEigvlsT} evaluated at
  $\lambda=0$, i.e., the nonzero eigenvalues of
  $E(\alpha_1, \dots, \alpha_m)$ have so far all been described by the
  formula~\eqref{eq:rho}.
\end{remark}

\PL{
  \subsection{Structural invariance of $E$}
}

Since $\mathcal{V}_\lambda$ and $\mathcal{V}_{\mathbf{z}}$ are
invariant subspaces of $E$, so is their direct sum\footnote{Because
  $T$ is diagonalizable, the eigenvectors of $T$ are linearly
  independent, as are the columns of the matrices
  $V_{\lambda}, V_{\mathbf{z}}$; hence the direct sum is well
  defined.},
\begin{equation}\label{eqn_SpctrlAnal_VasDirectSum}
  \mathcal{V} := 
  \underbrace{  \left[ \bigoplus\limits_{\lambda^2 \neq 0} \mathcal{V}_{\lambda} \right] \oplus 
    \left[ \bigoplus\limits_{\mathbf{z} \notin \mathrm{Ker}(C)} \mathcal{V}_{\mathbf{z}} \right]  }_{=:  \mathcal{V}_1} \oplus 
  \underbrace{  \left[ \bigoplus\limits_{\mathbf{z} \in \mathrm{Ker}(C)} \mathcal{V}_{\mathbf{z}} \right]  }_{=:  \mathcal{V}_2}.
\end{equation}
The subspace $\mathcal{V}$ captures the already established
eigenspaces from
Propositions~\ref{prop_SpctrlAnal_NonZeroEigvlsT},~\ref{prop_SpctrlAnal_ZeroEigvlsT}
\PL{and~\ref{prop_SpctrlAnal_ZeroEigvlsT_KerC}}. \PL{Next, we study
  $E(\alpha_1,\dotsc ,\alpha_m)$ on the complement of $\mathcal{V}$.}

We set $k_0 := \mathrm{dim}( \mathrm{Ker}(C) )$ and notice that
\begin{equation*}
  \Ker = \left\{
    \begin{bmatrix}
      -A^{-1} B \vec{x} \\
      \vec{x}
    \end{bmatrix} \bigg \vert  \vec{x} \in \mathbb{C}^n \right\},
\end{equation*}
and using these we obtain $\mathrm{dim}(\mathcal{V}) = 2(N-n) - k_0$,
since each eigenvector of $T\in \mathbb{C}^{(N-n)\times (N-n)}$
induces two different eigenvectors of $E(\alpha_1, \dots, \alpha_m)$,
except for those lying in $\mathrm{Ker}(C)$, which induce only one.

\begin{proposition}\label{prop_SpctrlAnal_InvSbspc}
  Having $\mathcal{V}_1,\mathcal{V}_2$ as
  in~\eqref{eqn_SpctrlAnal_VasDirectSum}, we have
  \begin{equation*}
    \PL{\mathcal{V}_2} \cap \Ker = \{\mathbf{0}\}
    \quad \mathrm{and} \quad
    \mathrm{dim}(\PL{\mathcal{V}_1} \cap \Ker) = N-n-k_0.
  \end{equation*}
  \PL{Moreover, we can find $\mathcal{W}\subset \Ker$ such that
  $\mathbb{C}^N = \mathcal{V} \oplus \mathcal{W}$.}
\end{proposition}

\begin{proof}
  Recall that any element in \PL{$\mathcal{V}_2$} is of the form
  \begin{align*}
    \begin{bmatrix}
      \vec{z} \\
      ~
    \end{bmatrix}, \quad \vec{z} \in \mathrm{Ker}(C)
  \end{align*}
  and so the product of this vector with $E^{(c)}$ is then itself,
  which is nonzero, yielding the first part of the statement.

  For the second part, as \PL{$\mathcal{V}_1$} is a direct sum of
  smaller subspaces, we can consider each subspace individually.
  First, consider a vector $\vec{x} \in \mathcal{V}_\lambda \cap \Ker$
  for some $\lambda^2$ a nonzero eigenvalue of $T$ with eigenvector
  $\vec{w}$.  By definition of $\mathcal{V}_{\lambda}$, $\vec{x}$ has
  the form
  \begin{equation*}
    \vec{x} =
    \begin{bmatrix}
      \zeta_1 \vec{w} \\
      \zeta_2 D^{-1} C \vec{w}
    \end{bmatrix},
  \end{equation*}
  for some scalars $\zeta_1, \zeta_2$.  Moreover, since
  $\vec{x} \in \Ker$,
  \begin{equation}\label{eqn_SpctrlAnal_zeta1w}
    \zeta_1 \vec{w} = -\zeta_2 A^{-1} B D^{-1} C \vec{w} = -\zeta_2 T
    \vec{w} = - \zeta_2 \lambda^2 \vec{w}.
  \end{equation}
  Thus, $\zeta_1 = -\lambda^2 \zeta_2$.  The subspace
  $\mathcal{V}_\lambda \cap \Ker$ is then a one-dimensional subspace.
  The same argument holds when replacing $\mathcal{V}_\lambda$ with
  $\mathcal{V}_{\mathbf{z}}$ for some
  $\vec{z} \notin \mathrm{Ker}(C)$, only in such \PL{a }case the
  right-hand side of~\eqref{eqn_SpctrlAnal_zeta1w} vanishes, implying
  $\zeta_1 = 0$ and leaving $\zeta_2$ arbitrary.  As a result, the
  subspace \PL{$\mathcal{V}_1$} $\cap \Ker$ is then the direct sum of
  these one-dimensional subspaces, of which there are $N-n-k_0$,
  yielding the second equality of the statement. \PL{The proof is
    finalized by using the Sum and Intersection Theorem for vector
    spaces, which gives}
  \begin{equation*}
    \mathrm{dim}(\mathcal{V} \cap \Ker) = N-n-k_0 \implies \mathrm{dim}(\mathcal{V} + \Ker) = N \PL{\implies \mathbb{C}^N = \mathcal{V} + \Ker.}
  \end{equation*}
\end{proof}

\PL{
  Last but not least, we show that the calculation leading to
  Proposition~\ref{prop_SpctrlAnal_NonZeroEigvlsT} generalizes to
  $\mathbb{C}^N$. This is not surprising, if we recall
  Remark~\ref{rmk: rho_lambda} that highlights the fact that
  Proposition~\ref{prop_SpctrlAnal_NonZeroEigvlsT} is in a sense
  representative of the general picture.  To do that we generalize the
  definition of the symbol $\rho_{\lambda^2}$ from a scalar parametric
  polynomial in $\lambda^2$ (see
  Remark~\ref{rmk_SpctrlAnal_rhoLambdaSqrd}) to a general parametric
  polynomial. In particular, we set
  \begin{equation*}
    \rho_{T} := \frac{1}{2} \left(I + \sqrt{T}\right)
    \prod_{i=1}^m \left( (1 - \alpha_i)I - \alpha_i \sqrt{T} \right)^2 +
    \frac{1}{2} \left(I - \sqrt{T}\right)
    \prod_{i=1}^m \left( (1 - \alpha_i)I + \alpha_i \sqrt{T} \right)^2
  \end{equation*}
  The following proposition establishes~\eqref{eq:rho} on a block-matrix
  level.
}

\begin{proposition}\label{prop_SpctrlAnal_Esqrd}
  \PL{If} $T$ is diagonalizable, then
  \PL{
    \begin{equation*}
      E^{(c)}
      E^{(s)} E^{(s)} E^{(c)} =
      \begin{bmatrix}
        \rho_T & \rho_T A^{-1}B \\
               &
      \end{bmatrix},
    \end{equation*}
    where $E^{(s)} = \prod_i E^{(s)}(\alpha_i)$.
  }
\end{proposition}
\PL{
  \begin{proof}
  First, a direct calculation gives
  \begin{equation*}
    E^{(c)} E^{(s)} E^{(s)} E^{(c)} = 
    \begin{bmatrix} I \\ 0 \end{bmatrix} 
    \underbrace{ \begin{bmatrix} I & A^{-1}B \end{bmatrix} E^{(s)} E^{(s)} \begin{bmatrix} I \\ 0 \end{bmatrix} }_{=: F}
    \begin{bmatrix} I & A^{-1}B \end{bmatrix},
  \end{equation*}
  with $F\in\mathbb{C}^{(N-n)\times (N-n)}$, i.e., it remains to show
  $F=\rho_T$. To that end, let $Q$ be the eigenmatrix of $T$ so that
  $T= Q\Lambda^2 Q^{-1}$; the columns of $Q$ are composed of the
  vectors $\vec{w}$ and $\vec{z}$, see
  Remark~\ref{rmk_secSpctrlAnal_Notation}. We start by calculating the
  action of $F$ in the basis $Q$:
  \begin{equation*}
    \begin{aligned}
      F\vec{w} &= 
                 \begin{bmatrix}
                   I & A^{-1}B
                 \end{bmatrix}
                 E^{(s)} E^{(s)}
                 \begin{bmatrix}
                   \vec{w} \\
                   \vec{0}
                 \end{bmatrix} = 
                 \begin{bmatrix}
                   I & A^{-1}B
                 \end{bmatrix}
                 E^{(s)} E^{(s)} \left(\frac{1}{2\lambda} V_{\lambda}
                 \begin{bmatrix}
                   1 \\
                   1
                 \end{bmatrix} \right) =
                 \frac{1}{2\lambda}
                 \begin{bmatrix}
                   I & A^{-1}B
                 \end{bmatrix}
                 V_{\lambda}
                 \begin{bmatrix}
                   \rho_-^2 & \\
                            & \rho_+^2
                 \end{bmatrix}
                 \begin{bmatrix}
                   1 \\
                   1
                 \end{bmatrix} \\
               &= \frac{1}{2\lambda}
                 \begin{bmatrix}
                   (\lambda + \lambda^2)\vec{w} & (\lambda - \lambda^2)\vec{w}
                 \end{bmatrix}
                 \begin{bmatrix}
                   \rho_-^2 & \\
                            & \rho_+^2
                 \end{bmatrix}
                 \begin{bmatrix}
                   1 \\
                   1
                 \end{bmatrix} =
                 \rho_{\PL{\lambda^2}} \vec{w}, \\
      F\vec{z} &= 
                 \begin{bmatrix}
                   I & A^{-1}B
                 \end{bmatrix}
                 E^{(s)} E^{(s)}
                 \begin{bmatrix}
                   \vec{z} \\
                   \vec{0}
                 \end{bmatrix} = 
                 \begin{bmatrix}
                   I & A^{-1}B
                 \end{bmatrix}
                 E^{(s)} E^{(s)} V_{\vec{z}}
                 \begin{bmatrix}
                   0 \\
                   1
                 \end{bmatrix} =
                 \begin{bmatrix}
                   I & A^{-1}B
                 \end{bmatrix}
                 V_{\vec{z}} \prod\limits_{i=1}^{m}
                 \begin{bmatrix}
                   (1-\alpha_i) & -\alpha_i \\
                                & (1-\alpha_i)
                 \end{bmatrix}^2
                 \begin{bmatrix}
                   0 \\
                   1
                 \end{bmatrix} \\
               &= \vec{z}
                 \begin{bmatrix}
                   0 & 1
                 \end{bmatrix} \left( \prod\limits_{i=1}^{m}
                 \begin{bmatrix}
                   1-\alpha_i & -\alpha_i \\
                              & 1-\alpha_i
                 \end{bmatrix}^2 \right)
                 \begin{bmatrix}
                   0 \\
                   1
                 \end{bmatrix}
                 = \prod\limits_{i=1}^{m} (1-\alpha_i)^2 \vec{z} = \rho_{\lambda^2=0} \vec{z}.
    \end{aligned} 
  \end{equation*}
  Hence,
  $F = Q~\mathrm{diag}\left(\rho_{\lambda_1^2},\dotsc ,
    \rho_{\lambda_{(N-n)}^2}\right)~Q^{-1} =
  Q~\rho_{\Lambda^2}~Q^{-1}$ and since $\rho_{\Lambda^2}$ is a
  polynomial, we obtain
  $Q~\rho_{\Lambda^2}~Q^{-1} = \rho_{Q~\Lambda^2~Q^{-1}} = \rho_{T}$,
  finishing the proof.
\end{proof}}

\PL{
  This is as far as we go without concrete choices of the smoothing
  parameters $\alpha_1,\dotsc ,\alpha_m$.  We have explicitly
  calculated $2(N-n) - k_0$ eigenpairs of $E$ and we let a convenient
  choice of $\alpha_1,\dotsc ,\alpha_m$ do the rest of the job for us,
  i.e., show that the remaining $2n + k_0 - N$ eigenpairs of $E$
  exist.
}

\section{Spectral clustering}\label{sec:spectralclustering}

In this section we compute the optimal smoothing parameters
$\hat{\alpha}_1,\dotsc ,\hat{\alpha}_m$.  \PL{Propositions
\ref{prop_SpctrlAnal_NonZeroEigvlsT}--\ref{prop_SpctrlAnal_ZeroEigvlsT_KerC}}
give the nonzero eigenvalues of the error propagation operator $E$ as
\begin{equation*}
  \rho_{\PL{\lambda^2}} = \frac{1}{2} (1 + \lambda) \prod_{i=1}^m (1 - \alpha_i - \alpha_i \lambda)^2 + \frac{1}{2} (1 - \lambda) \prod_{i=1}^m (1 - \alpha_i + \alpha_i \lambda)^2,
\end{equation*}
where $\lambda^2$ is a (possibly zero) eigenvalue of $T$.  If $m=1$,
the nonzero eigenvalues of $E$ become
\begin{equation*}
  \rho_{\PL{\lambda^2}}(\alpha_1) = (1-\alpha_1)^2 + \alpha_1 (3\alpha_1-2) \lambda^2  .
\end{equation*}
From the viewpoint of choosing a suitable $\alpha_1$,
$\rho_{\PL{\lambda^2}}(\alpha_1)$ is a parametric polynomial whose
values give the nonzero eigenvalues of $E$. If we were able to choose
$\alpha_1 \in \mathbb{C}$ so as to nullify
$\rho_{\PL{\lambda^2}}(\alpha_1)$ for all possible \PL{$\lambda^2$},
then $E$ would become the zero operator and hence $M^{-1} = L^{-1}$,
see~\eqref{eqn:fullerror} and below. Unfortunately, such a choice of
$\alpha_1$ clearly does not exist.  Instead, we seek to perfectly
cluster all the nonzero eigenvalues of $E$ to a single, nonzero value
so that $E$ has in total two distinct eigenvalues.  To do so, we make
$\rho_{\PL{\lambda^2}}(\alpha_1)$ independent of \PL{$\lambda^2$},
obtaining $\hat{\alpha}_1 =2/3$.

Similarly simple formulas can be found for general $m$; to do so, we
first present a convenient identity, based on~\citep[\S 1.39, 1.396(3)
and 1.396(2)]{GradshteynRyzhik2007}.

\begin{lemma}\label{lemma_secSpctrlClstrin_TrigIdentity}
  Let $m\ge 1$ and define $\theta_i := \frac{2\pi i}{2m+1}$ for
  $i=1,\dotsc,m$. Then, for any \PL{$x\in \mathbb{C}\setminus\{0,\pm1\}$},
  \begin{equation*}
    \prod_{i=1}^{m}\left( \cos\theta_i \pm \frac{1}{2} (x+x^{-1}) \right)
    = \PL{(\pm 1)^m} \frac{x^{2m+1}\pm 1}{2^{m} x^{m} (x\pm 1)}.
  \end{equation*}
\end{lemma}
\begin{proof}
  Set $t=\tfrac{1}{2}(x+x^{-1})$ for convenience.
  For each $i$,
  \begin{equation*}
    \cos\theta_i \pm t
    = \cos\theta_i \pm \frac{x+\frac{1}{x}}{2}
    = \frac{2 \cos\theta_i x \pm \left( x^2 + 1 \right)}{2x}
    = \PL{\pm}\frac{x^2 \pm 2 x \cos\theta_i + 1}{2x}.
  \end{equation*}
  Taking the product over $i=1,\dots,m$ gives
  \begin{equation*}
    \prod_{i=1}^{m}\left( \cos\theta_i \pm t \right)
    = \prod_{i=1}^{m} \PL{\pm}
      \frac{x^2 \pm 2 x \cos\theta_i + 1}{2x}
    = \frac{\PL{(\pm 1)^m}}{\left( 2x \right)^m}
      \prod_{i=1}^{m} \left( x^2 \pm 2 x \cos\theta_i + 1 \right).
  \end{equation*}
  By \citep[\S1.39, 1.396(3) and 1.396(2)]{GradshteynRyzhik2007}
  \begin{equation*}
    \prod_{i=1}^{m}\left( x^2 + 2 x \cos\theta_i + 1 \right) = \frac{x^{2m+1}+1}{x+1}, \quad
    \prod_{i=1}^{m}\left( x^2 - 2 x \cos\theta_i + 1 \right) = \frac{x^{2m+1}-1}{x-1}.
  \end{equation*}
  Substituting each case yields the claim.
\end{proof}
Equipped with Lemma~\ref{lemma_secSpctrlClstrin_TrigIdentity}, we can
formulate the key result of this section.

\begin{proposition}\label{prop_SpctrlClstrin_RhoClstrn}
  Let $m\ge1$ and $\rho_{\PL{\lambda^2}}$ be a polynomial as defined in
  equation \eqref{eq:rho}.  Let
  \begin{equation*}
    c_i =-\cos \left ( \frac{2\pi i}{2m+1} \right )
    \quad \PL{\mathrm{and} \quad
    \hat{\alpha}_i := \frac{1}{1+c_i}.}
  \end{equation*}
  Then
  \begin{equation*}
    \hat{\rho}_m   := \rho_{\PL{\lambda^2}}(\hat{\alpha}_1, \dots,
    \hat{\alpha}_m) = \frac{1}{(2m+1)^2}.
  \end{equation*}
\end{proposition}

\begin{proof}
  It suffices to show that the polynomial $\rho_{\PL{\lambda^2}}$ evaluated
  at this choice of $\alpha_i=1/(1+c_i)$ is equal to the stated
  formula, as it is then independent of $\lambda$.

  Rearrange the formula for $\hat{\alpha}_i$ to isolate for $c_i$:
  \begin{equation*}
    c_i = \frac{1-\hat{\alpha}_i}{\hat{\alpha}_i} \implies
    (1-\hat{\alpha}_i)\mp \hat{\alpha}_i  \lambda = \hat{\alpha}_i(c_i\mp  \lambda),
  \end{equation*}
  and we can write
  \begin{equation}\label{eqn_secSpctrlClstrin_PartOfProofRhomAsProd}
    \hat{\rho}_m =
    \left(\prod_{i=1}^{m}\hat{\alpha}_i\right)^{2} \left[\tfrac12(1+ \lambda)\prod_{i=1}^{m}(c_i- \lambda)^{2} + \tfrac12(1- \lambda)\prod_{i=1}^{m}(c_i+ \lambda)^{2}\right].
  \end{equation}
  Using the change of variables
  $ \lambda=\tfrac{1}{2}(x+x^{-1})$ with
  $x \in \mathbb{C}\backslash \PL{\{0,\pm 1\}}$ as in
  Lemma~\ref{lemma_secSpctrlClstrin_TrigIdentity} we have
  \begin{equation*}
    \prod_{i=1}^{m}(c_i- \lambda)^2 = \left( \frac{x^{2m+1}+1}{2^{m}x^{m}(x+1)} \right)^2 \text{ and }
    \prod_{i=1}^{m}(c_i+ \lambda)^2 = \left( \frac{x^{2m+1}-1}{2^{m}x^{m}(x-1)} \right)^2.
  \end{equation*}
  Substituting back into
  eq. \eqref{eqn_secSpctrlClstrin_PartOfProofRhomAsProd} eventually
  yields
  \begin{equation*}
    \begin{aligned}
      \hat{\rho}_m
      =& \left(\prod_{i=1}^{m}\hat{\alpha}_i\right)^{2} \left\{\tfrac{1}{2} \left(1+\tfrac{x+\frac{1}{x}}{2}\right) \left(\frac{x^{2m+1}+1}{2^{m}x^{m}(x+1)}\right)^{2} + 
         \tfrac{1}{2} \left(1-\tfrac{x+\frac{1}{x}}{2}\right) \left(\frac{x^{2m+1}-1}{2^{m}x^{m}(x-1)}\right)^{ 2} \right\} \\
      =& \left(\prod_{i=1}^{m}\hat{\alpha}_i \right)^{2} \left\{ \frac{(x^{2m+1}+1)^{2}-(x^{2m+1}-1)^{2}}{2^{2m+2}x^{2m+1}} \right\}\\
      =& \left(\prod_{i=1}^{m}\hat{\alpha}_i\right)^{2} \left\{ \frac{\left[(x^{2m+1}+1)-(x^{2m+1}-1)\right]
         \left[(x^{2m+1}+1)+(x^{2m+1}-1)\right]}{2^{2m+2}x^{2m+1}}\right\} \\
      =& \left(\prod_{i=1}^{m}\hat{\alpha}_i\right)^{2} \left\{\frac{1}{2^{2m+2}x^{2m+1}}\cdot 2 \cdot 2x^{2m+1} \right\} 
         = \frac{\left(\prod_{i=1}^{m} \hat{\alpha}_i\right)^{2}}{2^{2m}} = \left( \frac{2^{-m}}{\PL{\prod_{i=1}^{m} \left(\cos \theta_i - 1\right)}} \right)^{2}.    
    \end{aligned}
  \end{equation*}
  We want to use Lemma~\ref{lemma_secSpctrlClstrin_TrigIdentity} with
  $x=1$, but this is an indeterminate limit of type 0/0, so first we
  must eliminate the two factors of $x-1$:
  \begin{equation*}
    \frac{x^{2m+1}-1}{2^mx^m(x-1)} = \frac{(x-1)\left(\sum_{j=0}^{2m} x^{2m-j} \right)}{2^mx^m(x-1)} = \frac{\sum_{j=0}^{2m} x^{2m-j}}{2^mx^m}.
  \end{equation*}
  We can now prove the statement, substituting $x=1$ in the previous
  equation:
  \begin{equation*}
    \rho_{\PL{\lambda^2}} = \left( 2^m \prod_{i=1}^{m} (\cos \theta_i - 1) \right)^{-2}
    = \left (2^m \frac{\sum_{j=0}^{2m} 1^{2m-j}}{2^m 1^m} \right )^{-2}
    = \frac{1}{\left(2m+1\right)^2}.
  \end{equation*}
\end{proof}

\PL{
\begin{theorem} \label{theo:cluster}
  Let $T=A^{-1}BD^{-1}C$ be diagonalizable and denote by
  $\hat{M}^{-1}$ the action of
    Algorithm~\ref{alg:TwoLevelVCycleClustered} with the smoothing
    coefficients $\hat{\alpha}_1,\dotsc ,\hat{\alpha}_m$. Then the
  error propagation operator
  $E(\hat{\alpha}_1, \dots, \hat{\alpha}_m) \equiv \hat{E}$ is also
  diagonalizable and so is the preconditioned system $\hat{M}^{-1}L$.
  Moreover,
  \begin{equation*}
    \sigma \left( \hat{E} \right) = \left\{ 0, \frac{1}{(2m+1)^2} \right\}
    \quad \mathrm{and} \quad
    \sigma \left( \hat{M}^{-1}L \right) = \left\{ 1, 1-\frac{1}{(2m+1)^2} \right\},
  \end{equation*}
  with the eigenvectors of $\hat{E}$ corresponding to $1/(2m+1)^2$ given by
  Propositions~\ref{prop_SpctrlAnal_NonZeroEigvlsT},~\ref{prop_SpctrlAnal_ZeroEigvlsT}
  and~\ref{prop_SpctrlAnal_ZeroEigvlsT_KerC}. Thus, if these
  assumptions hold on every level, Algorithm~\ref{alg:AMGK} is direct
  in exact arithmetic.
\end{theorem}
\begin{proof}
  As a consequence of Proposition~\ref{prop_SpctrlClstrin_RhoClstrn},
  $\hat{\rho}_m \equiv \rho_{\lambda^2}(\hat{\alpha}_1,\dotsc ,
  \hat{\alpha}_m)$ is independent of $\lambda^2$. Hence, evaluating
  the polynomial
  $\rho_{\lambda^2}(\hat{\alpha}_1,\dotsc , \hat{\alpha}_m)$ at the
  matrix $T$, the resulting matrix will be also independent of
  $T$:
  \begin{equation*}
    \rho_{T}(\hat{\alpha}_1,\dotsc , \hat{\alpha}_m) = \hat{\rho}_m I,
  \end{equation*}
  which, thanks to Proposition~\ref{prop_SpctrlAnal_Esqrd}, gives the
  identity $\hat{E}^2 = \hat{\rho}_m \hat{E}$, or, equivalently,
  \begin{equation*}
    \hat{E}^2 - \hat{\rho}_m \hat{E} = \hat{E} \left( \hat{E} - \hat{\rho}_m I\right) = 0
    \implies \hat{p}(t) := t(t-\hat{\rho}_{m})
  \end{equation*}
  is the minimal polynomial of $\hat{E}$.  As $\hat{p}(t)$ has
  distinct roots, we conclude that $\hat{E}$ is diagonalizable.
  Moreover, the spectrum of $\hat{E}$ is given by the roots of
  $\hat{p}(t)$. Recalling that $\hat{M}^{-1}L = I-\hat{E}$, we obtain
  $\sigma(\hat{M}^{-1}L) = 1 - \sigma(\hat{E})$. Since
  $\mathrm{rank}(\hat{E}) \le \mathrm{rank}(E^{(c)}) = N-n$, the $N-n$
  linearly independent eigenvectors for $\hat{\rho}_m$ given by
  Propositions~\ref{prop_SpctrlAnal_NonZeroEigvlsT}--\ref{prop_SpctrlAnal_ZeroEigvlsT_KerC}
  span the corresponding eigenspace, and the last claim follows from
  Proposition~\ref{ref:KrylovConvergence} and Lemma~\ref{lem:AMGK}.
\end{proof}
}

\PL{
  We want to highlight the difference between the driving mechanism
  of Algorithm~\ref{alg:TwoLevelVCycle} on the one hand and
  Algorithm~\ref{alg:TwoLevelVCycleClustered} on the other.
}

\PL{
  The directness of Algorithm~\ref{alg:TwoLevelVCycle} comes from the
  fact that the error operator corresponding to the singular smoother
  $S^{-1}=\mathrm{blkdiag}(A^{-1},0)$ is a projection:
  \begin{equation*}
    S^{-1}LS^{-1}=S^{-1} \quad \text{and  so} \quad (I-S^{-1}L)(I-S^{-1}L) = I-S^{-1}L,
  \end{equation*}
  and the coarse space correction complements this projection to
  annihilate the error. In particular, repeating the singular smoother
  has no effect: it is not an iterative method. To make it one,
  approximations and/or damping have to be considered and it is not
  straightforward to extend the analysis to compensate for these.
}

\PL{
  Contrast that with Algorithm~\ref{alg:TwoLevelVCycleClustered}: the
  smoother corresponds to a standard stationary iterative method
  (which may or may not converge on its own). To compensate for the
  lack of the structural properties, we choose the smoothing
  coefficients carefully: instead of making the smoother converge on
  its own (as has been done before, see~\citep[Section
  5.5.2]{LiesenStrakos2013}), we choose them to complement the coarse
  correction. The resulting method is not direct but its error
  operator has a very enticing geometrical property:
  $\hat{E}^2 = \hat{\rho}_m \hat{E}$. In words, the origin, the error
  after one iteration and the error after two iterations must lie on a
  single line. Hence, after two iterations, we can extrapolate the
  error based on the approximations and compute the true
  solution. This is what the Krylov acceleration of
  Algorithm~\ref{alg:AMGK} (or the
  formula~\eqref{eqn_secDirMeth_LinvAsFunctnOfLMinv}) does.
}

\PL{
  Last, we want to demonstrate that Theorem~\ref{theo:cluster} does
  rely on the particular choice
  $\hat{\alpha}_1,\dotsc , \hat{\alpha}_m$ for the
  diagonalizability. For example, taking $A=D=1, B=C=3/\sqrt{5}, m=1$
  and $\alpha=1/4$, we obtain $T=9/5$ (diagonalizable) while $E$
  becomes similar to a Jordan cell:
  \begin{equation*}
    E = \frac{9}{200}
    \begin{bmatrix}
      5 & 5\sqrt{5} \\
      -\sqrt{5} & -5
    \end{bmatrix} = \frac{9}{200}
    \begin{bmatrix}
      -5\sqrt{5} & -\sqrt{5} \\
      5 & 0 \end{bmatrix}
    \begin{bmatrix}
      0 & 1 \\
      0 & 0
    \end{bmatrix}
    \begin{bmatrix}
      -5\sqrt{5} & -\sqrt{5} \\
      5 & 0
    \end{bmatrix}^{-1}.
  \end{equation*}
}

\section{Adaptation of the direct solver to practice}\label{sec:adaptation}
This section highlights some practical considerations and examples,
showcasing the potential of the above results beyond the abstract
setting.

\subsection{Conditioning of the preconditioned system}\label{sec_Appl_DirSolv_and_FoV}

We start by numerically illustrating Theorem \ref{theo:cluster}. As
our test piece, we construct a complex-valued, non-normal matrix
$L \in\mathbb{C}^{N\times N}$, for which we can directly control the
non-normality as well as \PL{the definiteness of its Hermitian
  part}. We do that by combining a \PL{Hermitian} matrix $H$ and a
skew-Hermitian matrix $K$, modulated by a parameter $\gamma >0$:
\begin{equation*}
L = H + \gamma K.
\end{equation*}
To construct $H$ and $K$, assemble $W,X\in \mathbb{C}^{N\times N}$ by
drawing the real and imaginary parts of the matrices from the
distribution $\mathcal{N}(0,1)$, independently for each entry.  Set
$H = W^\ast W + \eta I$ and $K = \tfrac{1}{2} (X - X^\ast)$. With high
probability, $H$ is HPD for $\eta > 0$ and Hermitian indefinite for
$\eta < 0$, and $K$ is skew-Hermitian.  Analogously, we also generate
a random, complex-valued right-hand side vector and solve the
resulting system with Algorithm~\ref{alg:AMGK}.

\begin{figure}[pos=t]
  \centering
  \includegraphics[width=.99\linewidth]{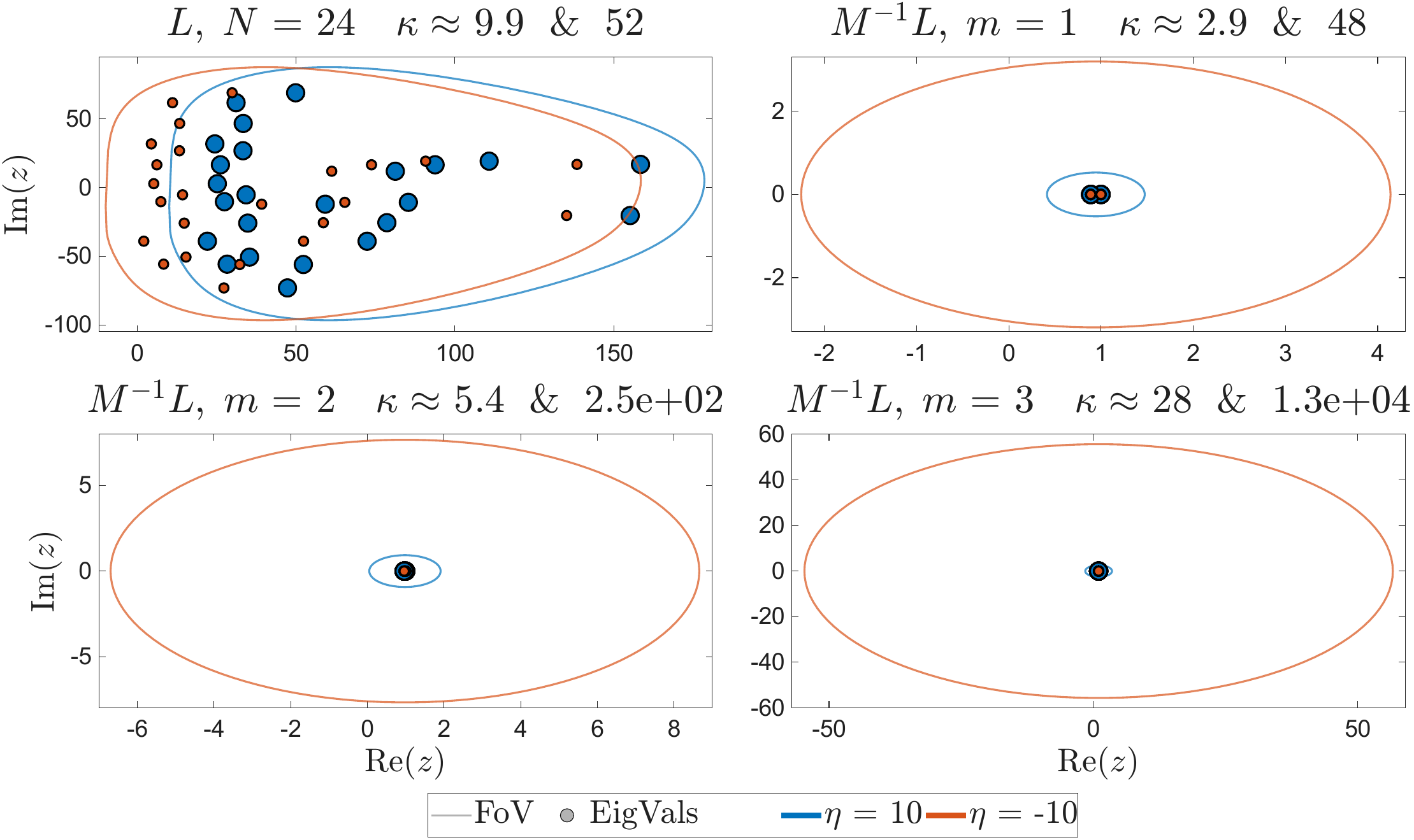}
  \caption{Field of values and spectra for $L$ (top-left) and
    \PL{$\hat M^{-1}L$} (rest) for $N=24$ with $\gamma=10$. We also give
    the condition number $\kappa$ of $L$ or \PL{$\hat M^{-1}L$} in the
    title of each of the subplots, first for $\eta=10$ and then
    for $\eta=-10$.}
  \label{fig_ApplFoV_EigvalsAndFoVBasedOnNmbSmthn}
\end{figure}

In Figure~\ref{fig_ApplFoV_EigvalsAndFoVBasedOnNmbSmthn},
we plot the spectra and the fields of values of $L$ and
\PL{$\hat M^{-1}L$} when \PL{the Hermitian part of $L$ is positive}
definite ($\eta=\gamma=10$) and indefinite ($\eta=-10, \gamma=10$),
and in all cases $L$ is well-conditioned ($\kappa \leq 10^2$). The
spectra of \PL{$\hat M^{-1}L$} are numerically clustered, specifically
on the level $n^{2m+1}\varepsilon_\text{mach}$, thus confirming the
results of the previous section.  Next, if \PL{the Hermitian part of
  $L$ is positive} definite, then the preconditioned system is also
well-conditioned even for larger numbers of pre- and post-smoothings.
However, when \PL{the Hermitian part of $L$ is} indefinite, the
preconditioned system becomes increasingly ill-conditioned as the
number of pre- and post-smoothings grows.  For $m=6$ we obtain
$\kappa\approx 10^8$ \PL{when the Hermitian part of $L$ is}
indefinite, compared to $\kappa\approx 10^4$ when \PL{the Hermitian
  part of $L$ is positive} definite.  This is dangerous, as
left-preconditioned GMRES measures the convergence by the relative
preconditioned residual.  For $m=3$ \PL{with an indefinite Hermitian
  part of $L$}, shown in the bottom-rightmost subplot, the relative
error of the resulting approximate solution becomes roughly of order
$10^{-10}$, while for \PL{a positive definite Hermitian part of $L$},
it stays at the level of $10^{-14}$.

We observe that the gap in conditioning of the preconditioned system
shrinks as $\eta\rightarrow 0^\pm$ and widens as
$\eta\rightarrow \pm\infty$.  The condition number grows steadily as
we increase $m$.  We note that from a practical viewpoint one is
usually not interested in running a larger number of pre- and
post-smoothings.  On the other hand, \PL{problems with an indefinite
  Hermitian part} are of practical importance and we advise caution as
we believe that further understanding is necessary.

\subsection{Discontinuous Galerkin 1D}\label{sec_Appl_DG}
We now look at applying the definition of the operators used in
Algorithm \ref{alg:TwoLevelVCycleClustered} to a discontinuous
Galerkin problem. Similar work to Section \ref{sec:spectralclustering}
has recently been explored from a geometric point of view
in~\citep{LuceroLorcaGander2022} for the discontinuous Galerkin method
applied to a simple test problem.  There, the authors find optimal
parameters within a family of prolongation and restriction operators
designed under strict locality and symmetry constraints.  In the
current algebraic setting, these geometric constraints are lifted.
Instead, we derive the prolongation operator directly from the
algebraic block structure of the reordered fine-grid matrix.

We perform \PL{a red-black permutation} so that the operator takes the
form
\begin{equation*}
  L =
  \begin{bmatrix}
    A & B \\
    C & D
  \end{bmatrix}, \quad
  P_{\mathrm{RB}} =
  \begin{bmatrix}
    - A^{-1}B \\
    I
  \end{bmatrix}, \quad
  R_{\mathrm{RB}} = P_{\mathrm{RB}}^{\top}, \quad
  M_0 = R_{\mathrm{RB}} L P_{\mathrm{RB}}.
\end{equation*}
We consider two admissible $1$D red-black reorderings:
\begin{itemize}
\item element-wise, grouping DoFs by alternating elements, and;
\item interface-wise, alternating traces across interfaces.
\end{itemize}
These lead to different local algebraic stencils for $P$ and for the
Jacobi smoother $S^{-1}$.  In all cases, the smoothing relaxation
parameters are $\hat{\alpha}_i$ defined in Proposition
\ref{prop_SpctrlClstrin_RhoClstrn}.

\begin{definition}[Element-wise local stencils]
  On each two-element patch, the local prolongation block (fine $8$
  DoFs from $4$ coarse DoFs) and the Jacobi smoother, which uses
  $2\times2$ diagonal blocks, are
  \begin{align}
    P_{\mathrm{loc,el}}(\delta)=
    \left(
    \begin{smallmatrix}
      \frac{\delta}{4\delta-2} & \frac{\delta-1}{4\delta-2} & 0 & 0\\ 
      \frac{\delta-1}{4\delta-2} & \frac{\delta}{4\delta-2} & 0 & 0\\ 
      1 & 0 & 0 & 0\\
      0 & 1 & 0 & 0\\ 
      \frac{\delta}{4\delta-2} & \frac{\delta-1}{4\delta-2} & \frac{\delta}{4\delta-2} & \frac{\delta-1}{4\delta-2}\\ 
      \frac{\delta-1}{4\delta-2} & \frac{\delta}{4\delta-2} & \frac{\delta-1}{4\delta-2} & \frac{\delta}{4\delta-2}\\ 
      0 & 0 & 1 & 0\\
      0 & 0 & 0 & 1
    \end{smallmatrix}
    \right),
    &\quad
      S^{-1}_{\mathrm{loc,el}}(\delta)
      =\frac{1}{-1+2 \delta}
      \left(
      \begin{smallmatrix}
        \delta & \delta-1\\
        \delta-1 & \delta
      \end{smallmatrix}
      \right),
  \end{align}
  \PL{where $\delta>1$ denotes the DG penalty parameter (analogous to
    $\delta_0$ in \citep{LuceroLorcaGander2022}) and assembly uses the
    usual overlap; the resulting global $P$ is \emph{not} block
    diagonal.}
\end{definition}

\begin{definition}[Interface-wise local stencils]
  On each interface, the local prolongation block and the Jacobi
  smoother are
  \begin{align}
    P_{\mathrm{loc,int}}(\delta)=
    \left(
    \begin{smallmatrix}
      \frac{\delta-1}{\delta} & \frac{1}{2\delta} & 0 & 0\\
      1 & 0 & 0 & 0\\
      0 & 1 & 0 & 0\\
      \frac{1}{2\delta} & \frac{\delta-1}{\delta} & \frac{1}{2\delta} & 0\\
      0 & \frac{1}{2\delta} & \frac{\delta-1}{\delta} & \frac{1}{2\delta}\\
      0 & 0 & 1 & 0\\
      0 & 0 & 0 & 1\\
      0 & 0 & \frac{1}{2\delta} & \frac{\delta-1}{\delta}
    \end{smallmatrix}
    \right),
    & \quad
      S^{-1}_{\mathrm{loc,int}}(\delta)=\frac{1}{\delta} I_{2},
  \end{align}
  \PL{where $\delta>1$ denotes the DG penalty parameter.}
\end{definition}

The definition of $R_{\mathrm{RB}}$ and $M_0$ follows directly from
the symmetric two-level formulation introduced in
Section~\ref{sec:Background}.  The theory developed in
Section~\ref{sec:spectralclustering} applies for an arbitrary number
of smoothing steps $m$.  Obtaining such a closed-form description of
the spectrum for general $m$ would be extremely cumbersome with a
classical LFA approach, which requires a mode-by-mode treatment and
becomes rapidly intractable as the number of smoothing steps
increases.

Compared with the interpolation from \citep{LuceroLorcaGander2022},
this algebraic prolongation has a richer coupling pattern: some
fine-grid degrees of freedom are interpolated from two neighboring
coarse elements instead of one.  This stronger coupling increases the
effective range of the operator and removes the need for parameter
tuning.  The two-point spectral clustering, which in
\citep{LuceroLorcaGander2022} appeared only for a specific triple
$(\alpha,c,\delta_0)$, now emerges automatically for any
\PL{admissible} choice of $\delta$.

From the DG perspective, this prolongation is again discontinuous.
The mixed rows of $P_{\mathrm{loc}}(\delta)$ generally produce
different fine values on the two sides of an interface.  A continuous
coarse function is typically mapped to a discontinuous fine one.
Together with the results of \citep{LuceroLorcaGander2022}, this
reinforces the same qualitative conclusion: When constructing
efficient multilevel preconditioners for DG discretizations, one
should use discontinuous interpolation operators, even when the
interpolation arises purely from the algebraic two-level factorization
rather than from an explicit geometric design.  This shows that the
clustering phenomenon does not depend on the geometric notion of
continuity or discontinuity of the interpolation, but follows directly
from the algebraic structure of the symmetric two-level formulation
itself.

\subsection{Sparsity of prolongation and restriction operators}\label{sec_Appl_FD}

Using Algorithm~\ref{alg:AMGK} with the choices of operators as in
eq.~\eqref{eqn:optimalchoice2}, the operators quickly become too
expensive to construct and apply. In this section we consider
adaptations that still outperform the widely used choices for a simple
model problem.

Consider the model problem of the two–dimensional Poisson equation
with homogeneous Dirichlet boundary conditions on the unit square,
discretizing $-\Delta$ on a uniform finite–difference mesh. The system
matrix $L\in \mathbb{R}^{N\times N}$ has a \emph{Kronecker sum}
structure, i.e., if $L_{1D}\in\mathbb{R}^{\sqrt{N}\times \sqrt{N}}$
denotes the one-dimensional second-order finite-difference Laplacian
corresponding to the stencil $\tfrac1{h^2}[-1~2~-1]$ on a uniform mesh
with $\sqrt{N}$ interior points and spacing $h$, we have
\begin{equation}\label{eqn_secApplFD_LeqL1DoplusL1D}
L = L_{1D} \oplus L_{1D} \equiv L_{1D} \otimes I + I \otimes L_{1D}.
\end{equation}
Applying Algorithm~\ref{alg:TwoLevelVCycleClustered} \emph{directly}
is possible but the prolongation and restriction operators then
require the construction of $-A^{-1}B$ (and its Hermitian conjugate),
which is a dense matrix, posing a significant bottleneck to the
efficiency of the multilevel method.

One common practice is to approximate expensive operators by \PL{what
  we call here} \emph{tensor lifting}: if the differential operator is
a tensor product or sum, then we can use this algebraically and
approximate the ideal operators by the tensor product or sum of their
lower-dimensional discretizations. \PL{It follows the view that, with
  tensor-product polynomials, the one-dimensional problem is the
  building block of the higher-dimensional
  case~\citep{HemkerHoffmannVanRaalte2004,
    HemkerHoffmannVanRaalte2003}, and one-dimensional results have
  been carried over to higher dimensions for interpolation
  operators~\citep{LuceroLorcaGander2022} and relaxation
  parameters~\citep{LuceroLorcaGander2024}.} We demonstrate this
heuristic on the model problem above. First, we red-black reorder the
DoFs of $L_{1D}$ so that
\begin{equation*}
  \Pi_{1D}^{\mathrm{RB}} L_{1D} (\Pi_{1D}^{\mathrm{RB}})^T =
  \begin{bmatrix}
    A_{1D} & B_{1D} \\
    B_{1D}^T & D_{1D}
  \end{bmatrix},
\end{equation*}
and $A_{1D},D_{1D} \in \mathbb{R}^{\sqrt{N}/2\times \sqrt{N}/2}$
become diagonal. Then, we apply the algebraic construction of
Algorithm~\ref{alg:TwoLevelVCycleClustered} to $L_{1D}$ based on its
above blocking and obtain $P_{1D}, R_{1D}, S_{1D}^{-1}$. The
approximate prolongation, restriction and smoothing operators for $L$
then follow as
\begin{equation}\label{eqn_secApplFD_PRSinv_as_TensorLiftFrm1D}
  P = P_{1D}\otimes P_{1D}, \quad R = R_{1D}\otimes R_{1D}
  \quad \mathrm{and} \quad S^{-1} = \PL{\frac{1}{4}} \left(S^{-1}_{1D}\oplus S^{-1}_{1D} \right).
\end{equation}
Clearly, the results of previous sections do not apply in that we do
not obtain a direct solver, but an iterative one.

\begin{figure}[pos=t]
  \centering
  \includegraphics[width=.85\linewidth]{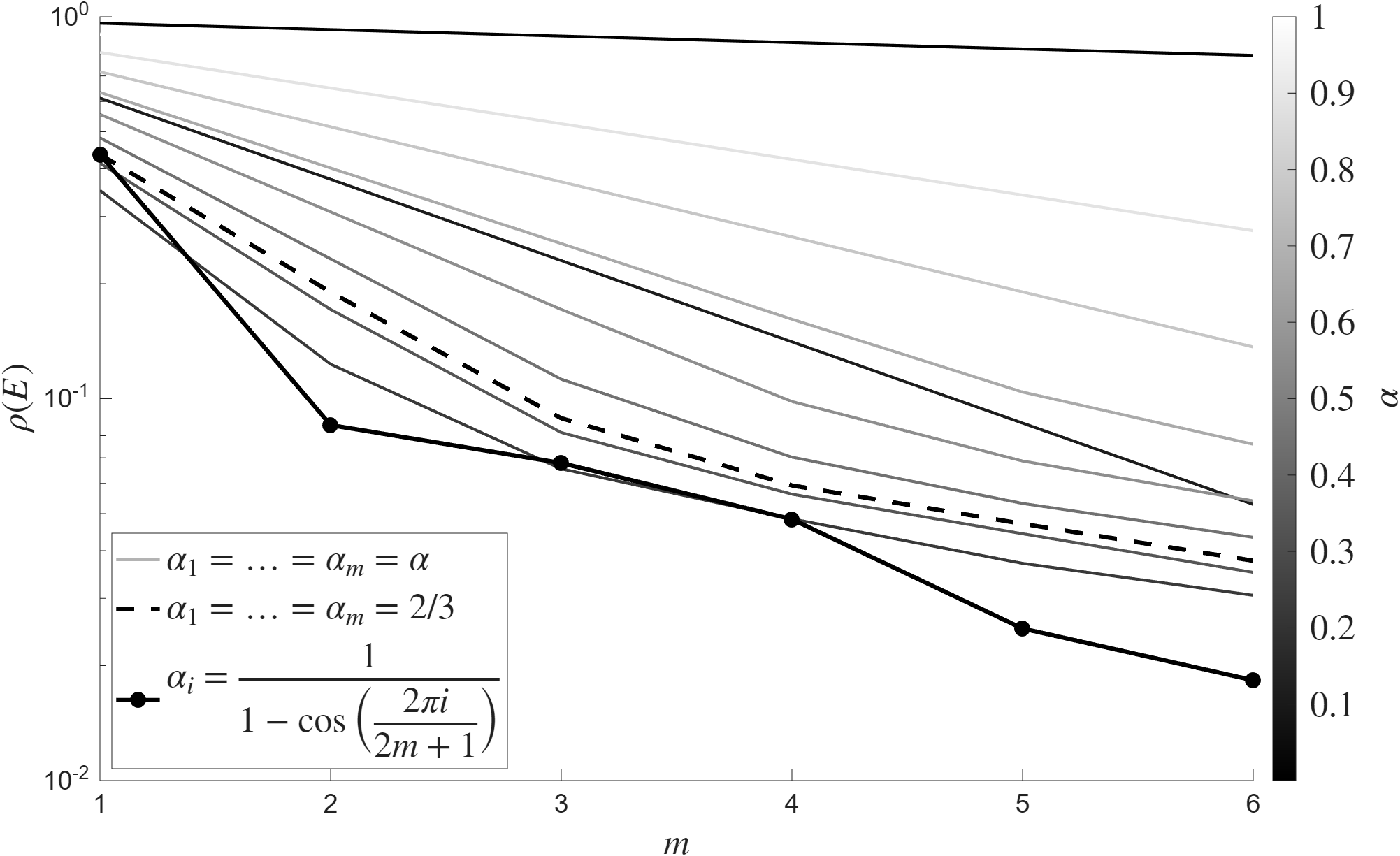}
  \caption{\PL{Spectral radius $\rho(E)$ as a function of the number
      of smoothing steps $m$ for the two-dimensional Laplacian system
      on a $16\times16$ interior finite-difference grid preconditioned
      with Algorithm~\ref{alg:TwoLevelVCycleClustered}. Gray curves:
      constant $\alpha\in\{0.1,\dots,1.0\}$; dashed: constant
      $\alpha=2/3$; thick curve with filled markers: smoothing
      parameters given by
      eq.~\eqref{eqn:optimalchoice3}.}} \label{fig_ApplFD_PoissonEqnSpectralRadius}
\end{figure}

The choice of the smoothing parameters $\alpha$ is usually done
heuristically, based on considering the smoother alone as
in~\citep{ErnstGander2013}, leading to the common choice
$\alpha_i = 2/3$ for all $i$, using LFA as in~\citep{Hackbusch1985} or
based on a lifting strategy similarly
to~\citep{LuceroLorcaGander2024}, where tedious analysis of the method
in 1D yields optimal parameters and these are simply used for
higher-dimensional problems. Adopting a similar mindset, we illustrate
in Figure~\ref{fig_ApplFD_PoissonEqnSpectralRadius} that even with the
choices in eq. \eqref{eqn_secApplFD_PRSinv_as_TensorLiftFrm1D} instead
of the ideal choices in eq. \eqref{eqn:optimalchoice2} our parameters
$\alpha_1,\dotsc ,\alpha_m$ perform \emph{near-optimally} compared to
other choices of $\alpha$. This suggests that the algebraic spectral
optimization of Section~\ref{sec:spectralclustering} extends
approximately to higher–dimensional tensor–product discretizations
while preserving sparse restriction and prolongation operators.

We also want to note that using either of these
restriction/prolongation operators, the coarse matrix $M_0$ becomes
denser than $L$ itself.  In our case, $M_0$ has the sparsity pattern
of the 9-point stencil and so neither of these techniques comes
without drawbacks in terms of efficiency of the multilevel method.

\section{Conclusion and future work}\label{sec:conclusion}

Algorithm \ref{alg:TwoLevelVCycleClustered} is a new algebraic
multigrid method with a clustered spectrum. Algorithm \ref{alg:AMGK},
its corresponding K-cycle, is then direct if used as a preconditioner
for an optimal Krylov method (or when used in the sense
of~\eqref{eqn_secDirMeth_LinvAsFunctnOfLMinv}) in exact arithmetic
after 2 iterations. The classic $2\times 2$-block partitioning inverse
is the basis for many methods in literature, e.g., AMG methods
sprouting from fine/coarse DoF grouping with approximations to the
ideal prolongations, restrictions and coarse spaces, Hierarchical
Poincaré--Steklov and other nested dissection methods. The algorithms
provide a best-of-class version of symmetric V- and W-cycle methods
and underline the usefulness of the K-cycle. The provided formulas for
the smoothing parameters function as a starting point in the design of
AMG methods.

\subsection{Krylov subspace methods}

Algorithm~\ref{alg:AMGK}, used as a preconditioner (inverting $M_0$
explicitly at the coarsest level), is a recursive, direct multilevel
method with the form of a K-cycle with two recursive coarse
corrections, thus also a W-cycle. Both the recursive application of
Algorithm~\ref{alg:TwoLevelVCycle} with the optimal choices in
eq. \eqref{eqn:optimalchoice1} and Algorithm \ref{alg:AMGK} produce
the exact solution after one \PL{and two iterations, respectively,}
when used as a preconditioner of a \PL{$\mathrm{KSM}^{\star}$ method},
but the former is a simpler V-cycle -- the K/W-cycle is the price to
pay for the ``symmetrization'' of the method and regularity of the
smoother.

Using CG/GMRES to precondition CG/GMRES has been considered before and
led to the flexible GMRES method (fGMRES)~\citep{Saad1993} and
flexible or generalized CG method~\citep{AxelssonVassilevski1991,
  Notay2000} (fCG, gCG). The term \emph{flexible} has been used to
highlight that the \PL{$\mathrm{KSM}^{\star}$ method} of choice does
not require the same preconditioner at each iteration; the
preconditioner can be nonlinear.

This becomes useful for numerical modifications or approximations of
(parts of) Algorithm~\ref{alg:AMGK}. As a result of these
approximations, two iterations of CG/GMRES no longer correspond to the
exact solver. Since Krylov subspace methods are generally nonlinear
solvers, we naturally obtain a nonlinear preconditioner within the
recursion. Therefore, we need to use the flexible versions of CG/GMRES
if we want to use an approximation of the multilevel version of
Algorithm~\ref{alg:AMGK} coupled with CG/GMRES.

\subsection{\PL{Recursive cross-separator elimination}}\label{sec_Concl_Crosses}

\PL{
  The singular choice~\eqref{eqn:optimalchoice1} in
  Algorithm~\ref{alg:TwoLevelVCycle} corresponds to the classic block
  factorization that is the algebraic backbone of
  substructuring~\citep[Section~6.1]{DoleanJolivetNataf2015}, nested
  dissection and HPS solvers~\citep{Martinsson2015, LuceroLorca2025,
    OutrataLuceroLorca2025}. This factorization is used both in direct
  solvers and within iterative methods, including on unstructured
  meshes~\citep{SchmitzYing2012,
    WildeyMuralikrishnanSriramkrishnanBuiThanh2019, LuLee2024} and for
  problems with non-separable coefficients~\citep{HoYing2016,
    SchmitzYing2012, MuralikrishnanBuiThanhShadid2020,
    WildeyMuralikrishnanSriramkrishnanBuiThanh2019}.  Under the
  assumptions of Theorem~\ref{theo:cluster}, Algorithm~\ref{alg:AMGK}
  inherits the requirements of Algorithm~\ref{alg:TwoLevelVCycle} and
  adds two conditions: the discretization must have no cross points,
  and applying $D^{-1}$ must cost no more than applying $A^{-1}$. The
  first condition holds for hybridized discontinuous Galerkin
  traces~\citep{MuralikrishnanBuiThanhShadid2020,
    WildeyMuralikrishnanSriramkrishnanBuiThanh2019}, cell-centered
  finite volumes~\citep{ChipmanCalhounBurstedde2024}, and spectral
  elements with the corner nodes
  removed~\citep{LuceroLorcaBeamsBeecroftGillman2024,
    LuceroLorca2025}. In all of these, after static condensation,
  every unknown lies on exactly one face, and two faces are coupled
  only if they belong to a common subdomain. The second condition is
  satisfied by recursive cross-separator elimination, in which the
  faces are grouped into red and blue crosses, each marked by a halo
  in Figure~\ref{fig_Appl_Crosses}. Since no two crosses of the same
  color share a subdomain, both $A$ and $D$ are block diagonal;
  moreover, no blue cross has more faces than a red one, so applying
  $D^{-1}$ costs no more than applying $A^{-1}$. Eliminating the red
  crosses leaves only the blue faces, which are the faces of the
  merged subdomains, and the construction is repeated on the next
  level. The red crosses correspond to the separator fronts
  in~\citep[Fig.~1]{MuralikrishnanBuiThanhShadid2020}
  and~\citep[Section~2.1]{SchmitzYing2012}, and to the interfaces
  eliminated by the 4-to-1 merges
  in~\citep{ChipmanCalhounBurstedde2024}.  On unstructured meshes, the
  subdomains can be obtained by overlaying a Cartesian grid, as
  in~\citep[Section~4.1]{SchmitzYing2012}.
}

\begin{figure}[pos=t]
  \centering
  \begin{tikzpicture}
    \def\pw{3}\def\gap{0.045}\def\nf{8}\def\hd{0.056}
    % \drawfaces{n}: the faces of the \nf x \nf decomposition that remain in
    % the n x n decomposition, with a gap at every cross point. A: faces on
    % the lines x,y = 1,3,... of the n x n decomposition, in crosses around
    % its (odd,odd) points; D: faces on its interior lines x,y = 2,4,..., in
    % crosses around its (even,even) points, truncated at the boundary.
    \def\drawfaces#1{%
      \pgfmathsetmacro\h{\pw/#1}\pgfmathsetmacro\hf{\pw/\nf}%
      \pgfmathtruncatemacro\nl{#1-1}%
      % halos: one pale stroke per cross, through its center and short of
      % its tips; red crosses are centered at (odd,odd) points, blue ones at
      % (even,even) points, truncated at the boundary.
      \foreach \ci in {0,...,#1}{\foreach \cj in {0,...,#1}{%
        \pgfmathtruncatemacro\pa{mod(\ci,2)}\pgfmathtruncatemacro\pb{mod(\cj,2)}%
        \ifnum\pa=\pb\relax
          \ifnum\pa=1\relax\def\hcol{red!25}\else\def\hcol{blue!25}\fi
          \pgfmathsetmacro\xa{max(0,(\ci-1)*\h)+\hd}%
          \pgfmathsetmacro\xb{min(\pw,(\ci+1)*\h)-\hd}%
          \pgfmathsetmacro\ya{max(0,(\cj-1)*\h)+\hd}%
          \pgfmathsetmacro\yb{min(\pw,(\cj+1)*\h)-\hd}%
          \ifnum\cj>0\relax\ifnum\cj<#1\relax
            \draw[color=\hcol,line width=3pt] (\xa,{\cj*\h}) -- (\xb,{\cj*\h});%
          \fi\fi
          \ifnum\ci>0\relax\ifnum\ci<#1\relax
            \draw[color=\hcol,line width=3pt] ({\ci*\h},\ya) -- ({\ci*\h},\yb);%
          \fi\fi
        \fi}}%
      \foreach \x in {1,...,\nl}{%
        \ifodd\x\relax\def\col{red}\else\def\col{blue}\fi
        \foreach \y in {1,...,\nf}{%
          \draw[color=\col,line width=1pt]
            ({(\y-1)*\hf+\gap},{\x*\h}) -- ({\y*\hf-\gap},{\x*\h})
            ({\x*\h},{(\y-1)*\hf+\gap}) -- ({\x*\h},{\y*\hf-\gap});}}%
      \foreach \y in {1,...,\nf}{%
        \draw[thick,gray!50]
          ({(\y-1)*\hf+\gap},0) -- ({\y*\hf-\gap},0)
          ({(\y-1)*\hf+\gap},\pw) -- ({\y*\hf-\gap},\pw)
          (0,{(\y-1)*\hf+\gap}) -- (0,{\y*\hf-\gap})
          (\pw,{(\y-1)*\hf+\gap}) -- (\pw,{\y*\hf-\gap});}}
    % Faces of the n x n decomposition (m = n/2) are numbered as follows: A
    % cross by cross, lexicographically in the center, arms E,N,W,S; then the
    % crosses of D likewise; then the single boundary faces of D (bottom,
    % left, right, top). \hface{X}{Y}: face on y=Y between x=X-1 and x=X;
    % \vface{X}{Y}: face on x=X between y=Y-1 and y=Y.
    \def\hface#1#2{%
      \ifodd#2\relax
        \pgfmathtruncatemacro\idx{4*((#2-1)/2*\m+ceil(#1/2)-1)+2*mod(#1,2)}%
      \else
        \pgfmathtruncatemacro\p{floor(#1/2)}\pgfmathtruncatemacro\q{#2/2}%
        \ifnum\p=0\relax\pgfmathtruncatemacro\idx{\nAD+\m-1+\q-1}%
        \else\ifnum\p=\m\relax\pgfmathtruncatemacro\idx{\nAD+2*(\m-1)+\q-1}%
        \else\pgfmathtruncatemacro\idx{\nA+4*((\q-1)*(\m-1)+\p-1)+2-2*mod(#1,2)}%
        \fi\fi
      \fi}
    \def\vface#1#2{%
      \ifodd#1\relax
        \pgfmathtruncatemacro\idx{4*((ceil(#2/2)-1)*\m+(#1-1)/2)+1+2*mod(#2,2)}%
      \else
        \pgfmathtruncatemacro\p{#1/2}\pgfmathtruncatemacro\q{floor(#2/2)}%
        \ifnum\q=0\relax\pgfmathtruncatemacro\idx{\nAD+\p-1}%
        \else\ifnum\q=\m\relax\pgfmathtruncatemacro\idx{\nAD+3*(\m-1)+\p-1}%
        \else\pgfmathtruncatemacro\idx{\nA+4*((\q-1)*(\m-1)+\p-1)+3-2*mod(#2,2)}%
        \fi\fi
      \fi}
    % \drawpattern{n}: one block row and column per face of the n x n
    % decomposition, each made of \nf/n faces of the \nf x \nf one; a block is
    % drawn when its two faces belong to a common subdomain.
    \def\drawpattern#1{%
      \pgfmathtruncatemacro\m{#1/2}\pgfmathtruncatemacro\nA{4*\m*\m}%
      \pgfmathtruncatemacro\nAD{\nA+4*(\m-1)*(\m-1)}%
      \pgfmathtruncatemacro\N{2*#1*(#1-1)}\pgfmathsetmacro\e{\pw/\N}%
      \begin{scope}[shift={(0,\pw)},scale=\e]
        \foreach \x in {1,...,#1}{\foreach \y in {1,...,#1}{%
          \def\fb{-1}\def\ft{-1}\def\fl{-1}\def\fr{-1}%
          \ifnum\y>1\relax\pgfmathtruncatemacro\z{\y-1}\hface{\x}{\z}\let\fb\idx\fi
          \ifnum\y<#1\relax\hface{\x}{\y}\let\ft\idx\fi
          \ifnum\x>1\relax\pgfmathtruncatemacro\z{\x-1}\vface{\z}{\y}\let\fl\idx\fi
          \ifnum\x<#1\relax\vface{\x}{\y}\let\fr\idx\fi
          \foreach \r in {\fb,\ft,\fl,\fr}{\foreach \c in {\fb,\ft,\fl,\fr}{%
            \ifnum\r>-1\relax\ifnum\c>-1\relax
              \ifnum\r<\nA\relax\ifnum\c<\nA\relax\def\col{red}\else\def\col{black}\fi
              \else\ifnum\c<\nA\relax\def\col{black}\else\def\col{blue}\fi\fi
              \fill[color=\col] (\c,-\r) rectangle ++(1,-1);
            \fi\fi}}}}
        \draw[gray] (0,0) rectangle (\N,-\N);
        \ifnum\nA<\N\relax\draw[gray] (\nA,0) -- (\nA,-\N) (0,-\nA) -- (\N,-\nA);\fi
      \end{scope}}
    \foreach \lev/\nsub in {1/8,2/4,3/2}{%
      \begin{scope}[shift={({(\lev-1)*(\pw+0.45)},0)}]
        \node[font=\footnotesize] at ({\pw/2},{2*\pw+0.55}) {Level \lev};
        \begin{scope}[shift={(0,\pw+0.3)}]\drawfaces{\nsub}\end{scope}
        \drawpattern{\nsub}
      \end{scope}}
  \end{tikzpicture}
  \caption{\PL{Top: faces remaining at levels 1--3 of an $8\times8$
    decomposition. Gaps mark the cross points, which carry no unknowns,
    and a pale halo groups the faces of each cross; red crosses form $A$
    and blue crosses form $D$, which is empty at level~3. Bottom: block
    sparsity pattern of the condensed operator at each level, ordered
    cross by cross with $A$ first; red and blue blocks belong to $A$ and
    $D$, respectively, and black blocks to $B$ and $C$.}}
  \label{fig_Appl_Crosses}
\end{figure}

\PL{
  As a direct solver, Algorithm~\ref{alg:AMGK} is more expensive than
  Algorithm~\ref{alg:TwoLevelVCycle}, since it additionally requires
  $D^{-1}$ and its K-cycle performs two recursive coarse corrections
  per level. Both methods share the main bottleneck of direct solvers,
  namely the dense Schur
  complement~\citep[Section~6.1]{DoleanJolivetNataf2015}, which
  motivates the use of approximations.
}

\PL{
  As mentioned below Theorem~\ref{theo:cluster},
  Algorithm~\ref{alg:TwoLevelVCycle} works based on a projection
  argument. In practice, approximate methods supplement or replace it
  with regular smoothers, such as block Jacobi before each level of
  skeletonization~\citep{FeliuFabaHoYing2020}, block Jacobi on the
  finest level~\citep{MuralikrishnanBuiThanhShadid2020}, and
  Jacobi-type, symmetric Gauss--Seidel, or block Jacobi smoothing on
  all levels in addition to the singular
  smoother~\citep{WildeyMuralikrishnanSriramkrishnanBuiThanh2019}. Whether
  smoothing on the coarser levels helps beyond elliptic problems is
  left open in~\citep{MuralikrishnanBuiThanhShadid2020}. We regard
  Algorithm~\ref{alg:AMGK} as an exact reference for methods with a
  regular smoother on every level. Its ingredients can be approximated
  as in the methods cited above, for instance by compressing the
  fronts~\citep{HoYing2016, SchmitzYing2012, LuLee2024}, by coarsening
  the mesh or projecting the traces onto coarser
  spaces~\citep{MuralikrishnanBuiThanhShadid2020,
    WildeyMuralikrishnanSriramkrishnanBuiThanh2019}, or by adjusting
  the smoothing parameters. The effect of such approximations on the
  two-point spectrum, and on the resulting methods, exceeds the scope
  of this manuscript and will be addressed in forthcoming
  publications.
}

\subsection{Newton--Schulz method}

An interesting way to interpret
eq. \eqref{eqn_secDirMeth_LinvAsFunctnOfLMinv} is to realize that the
damped Newton method for finding
$\left( \PL{\hat M^{-1}}L \right)^{-1}$ as the root of the function
$F(X) := X^{-1} - \PL{\hat M^{-1}}L$ reads
\begin{equation*}
  X_{i+1} = X_i + \omega \left(X_i - X_i \PL{\hat M^{-1}}LX_i\right),
\end{equation*}
with a damping parameter $\omega$. Taking $X_0=I$ as the initial guess
and the damping parameter as \PL{$\omega=1/(1-\hat{\rho}_m)$}, the
first iteration corresponds to
eq. \eqref{eqn_secDirMeth_LinvAsFunctnOfLMinv}. In other words,
eq. \eqref{eqn_secDirMeth_LinvAsFunctnOfLMinv} corresponds to one
iteration of the \PL{Newton--Schulz method with an optimal damping
factor $\omega$}. We obtain an iterative formulation that allows for
generalizations of the multilevel method in the case of approximations
or modifications. Moreover, the Newton link might allow us to tackle
some nonlinear operators either by finding nonlinear surrogates of the
ingredients of Algorithm \ref{alg:AMGK}, or by using it to approximate
the inverse of a Jacobian. Most approaches hinge on our capacity to
reconstruct the ingredients for a new iterate of the nonlinear method
-- we plan to expand on this in future work.

\PL{
\section*{Acknowledgments}
}

\PL{
  We would like to thank the two anonymous reviewers whose comments
  and remarks helped to meaningfully improve and correct the
  manuscript.
}

\bibliographystyle{cas-model2-names} 
\bibliography{paper} 

\end{document}